\documentclass[12pt,reqno,a4paper]{article}

\usepackage[blocks]{authblk}
\usepackage{amsmath,amsthm,amstext,amscd,amssymb,euscript,url}
\usepackage{mathrsfs,euscript,dsfont,calrsfs}
\usepackage{ushort}
\usepackage{comment}
\usepackage{relsize}
\usepackage{scalerel}
\usepackage{epsfig}
\usepackage{enumitem}
\usepackage{mathtools}
\usepackage{graphicx}
\usepackage{wrapfig}
\usepackage{caption}
\usepackage{multicol}
\usepackage{geometry}
\usepackage{hyperref}
\usepackage{todonotes}
\usepackage{biblatex}
\usepackage{comment}
\bibliography{refs}

\usepackage{xcolor}
\definecolor{unbleu}{rgb}{0.03, 0.15, 0.4}
\definecolor{monvert}{rgb}{0.0,.5,0.0}
\definecolor{britishracinggreen}{rgb}{0.0, 0.26, 0.15}
\definecolor{monbleu}{rgb}{0,.2,.8}
\definecolor{monautrebleu}{rgb}{0,0.4,.75}
\definecolor{applegreen}{rgb}{0.55, 0.71, 0.0}
\definecolor{monrouge}{rgb}{0.8, 0.0, 0.0} 
\definecolor{cadmiumgreen}{rgb}{0.0, 0.42, 0.24}
\definecolor{royalblue(traditional)}{rgb}{0.0, 0.14, 0.4}
\definecolor{black}{rgb}{0.0, 0.0, 0.0}
\definecolor{sepia}{rgb}{0.44, 0.26, 0.08}
\definecolor{teagreen}{rgb}{0.82, 0.94, 0.75}
\definecolor{yellow-green}{rgb}{0.6, 0.8, 0.2}
\definecolor{azure(colorwheel)}{rgb}{0.0, 0.5, 1.0}
\definecolor{awesome}{rgb}{1.0, 0.13, 0.32}
\definecolor{cadmiumyellow}{rgb}{1.0, 0.96, 0.0}
\definecolor{carrotorange}{rgb}{0.93, 0.57, 0.13}
\definecolor{green-yellow}{rgb}{0.68, 1.0, 0.18}
\definecolor{huntergreen}{rgb}{0.21, 0.37, 0.23}

\hypersetup{
pdfborder = {0 0 0},
colorlinks,
linkcolor=cadmiumgreen,
citecolor=cadmiumgreen,
urlcolor=cadmiumgreen
}

\makeatletter
\renewcommand{\tagform@}[1]{\maketag@@@{\color{cadmiumgreen}(#1)}}
\makeatother

\newcommand{\Z}{\mathds Z}
\newcommand{\R}{\mathds R}
\newcommand{\N}{\mathds N}

\newcommand{\E}{\mathds E}

\renewcommand{\phi}{\varphi}

\newcommand{\si}{\ensuremath{\sigma}}

\def\1{{\mathds 1}}

\newtheorem{theorem}{{\small T}{\scriptsize HEOREM}}[section]

\newtheorem{corollary}{{\bf{\small C}{\scriptsize OROLLARY}}}[section]
\newtheorem{proposition}{{\bf{\small P}{\scriptsize ROPOSITION}}}[section]
\newtheorem{lemma}{{\bf{\small L}{\scriptsize EMMA}}}[section]
\newtheorem{remark}{{\bf{\small R}{\scriptsize EMARK}}}[section]
\newtheorem{definition}{{\bf{\small D}{\scriptsize EFINITION}}}[section]

\newtheorem{induction}{{\bf{\small I}{\scriptsize NDUCTIVE HYPOTHESIS}}}[section]

\renewenvironment{proof}[1]
{\noindent{{\bf{\small{ P}{\scriptsize ROOF}}}.}\hspace{0.1cm} #1} {$\;\qed$\newline}

\newcommand{\be}{\begin{equation}}
\newcommand{\ee}{\end{equation}}

\newcommand{\bl}{\begin{lemma}}
\newcommand{\el}{\end{lemma}}

\newcommand{\br}{\begin{remark}}
\newcommand{\er}{\end{remark}}

\newcommand{\bt}{\begin{theorem}}
\newcommand{\et}{\end{theorem}}

\newcommand{\bd}{\begin{definition}}
\newcommand{\ed}{\end{definition}}

\newcommand{\bind}{\begin{induction}}
\newcommand{\eind}{\end{induction}}

\newcommand{\bp}{\begin{proposition}}
\newcommand{\ep}{\end{proposition}}

\newcommand{\bc}{\begin{corollary}}
\newcommand{\ec}{\end{corollary}}

\newcommand{\bpr}{\begin{proof}}
\newcommand{\epr}{\end{proof}}

\newcommand{\bi}{\begin{itemize}}
\newcommand{\ei}{\end{itemize}}

\newcommand{\ben}{\begin{enumerate}}
\newcommand{\een}{\end{enumerate}}

\newcommand{\caC}{{\mathscr C}}

\newcommand{\caN}{{\mathcal N}}

\newcommand{\p}{p^{\mathrm{rev}}}
\newcommand{\U}{U^{\mathrm{rev}}}
\newcommand{\lrev}{L_t^{\mathrm{rev}}}

\newcommand*\diff{\mathop{}\!\mathrm{d}}
\newcommand{\dert}{\frac{\mathop{}\!\mathrm{d}}{\mathop{}\!\mathrm{d} t}}
\newcommand{\ders}{\frac{\mathop{}\!\mathrm{d}}{\mathop{}\!\mathrm{d} s}}
\newcommand{\derx}{\frac{\mathop{}\!\mathrm{d}}{\mathop{}\!\mathrm{d} x}}
\newcommand{\pdert}{\frac{\mathop{}\!\partial}{\mathop{}\!\partial t}}
\newcommand{\pders}{\frac{\mathop{}\!\partial}{\mathop{}\!\partial s}}

\newcommand{\RW}{\mathrm{RW}}
\newcommand{\rev}{\mathrm{rev}}

\DeclareMathOperator{\e}{\mathrm{e}}

\begin{document}

\title{{\bf Duality and hidden symmetries for time-inhomogeneous Markov processes}}

\author[1,2]{Frank Redig
}
\author[1,2]{Josephine Schikhof
}

\affil[1]{Institute of Applied Mathematics, Delft University of Technology, Delft,
	The Netherlands}
\affil[2]{Dutch Institute for Emergent Phenomena, Netherlands}

\date{Dated: \today}

\maketitle

\begin{abstract}
	We study the notion of duality in the setting of time-inhomogeneous Markov processes.
	We define duality with time-dependent duality functions.
	We establish the connection between duality and time-dependent symmetries and intertwinings.
	We illustrate the abstract results with examples from interacting particle systems and diffusion processes.
\end{abstract}
\tableofcontents

\section{Introduction} Duality is a powerful tool in the study of interacting particle systems.
It enables to connect a system under study to a simpler dual system via a duality function.
A simple but useful illustration of duality is the fact that the expected number of particles of an interacting system such as the symmetric exclusion process (SEP) can be computed from the initial configuration and a single random walk.
The fact that this random walk scales to Brownian motion upon diffusive rescaling of space and time then gives for such systems the hydrodynamic limit in its simplest form, namely, in expectation.

In recent years, following pioneering work by Sandow and Sch\"{u}tz \cite{schutz1994non}, it became clear that systems which exhibit duality in some form often have special symmetries, i.e., operators commuting with the Markov generator (see \cite{giardina2009duality}, \cite{redig2026duality}).
More precisely, symmetries of the generator allow us to produce new duality functions from a given one.
Often a ``cheap'' duality function is available, e.g., via detailed balance.
This formalism, extensively studied in \cite{redig2026duality}, uses quite strongly that (up to domain issues) the Markov generator commutes with the semigroup, and therefore the construction of duality functions reduces to the study of the symmetries of the Markov generator.

In this paper, we extend these ideas to the time-inhomogeneous setting.
One of our motivations is the study of particle systems in a time-dependent geometry, where the particles are moving on a set of vertices which in turn move according to another interacting particle system.
A relevant setting is a neighborhood graph used to approximate a Riemannian manifold \cite{hein_graph_2007, van2020hydrodynamic}, where additionally one allows the vertices of this graph to move according to independent diffusion processes.
One can then study e.g., a symmetric exclusion process on these moving vertices and the hydrodynamic limit of the joint system of moving vertices and particles.
If one then fixes a realization of the moving vertices, the resulting process of the particles becomes a particle system in a dynamic environment, which is then unavoidably time-inhomogeneous.
Other natural examples are particle systems in dynamic random environment, where upon fixing (quenching) the evolution of the random environment, the resulting particle system is time-inhomogeneous (see e.g., \cite{redig2020dynamicSEP},\cite{seiler2023contact}).
If one is interested in duality properties of such systems, then one is led naturally into the question of extending the existing theory of Markov process dualities into the realm of time-inhomogeneous processes.

In this extension problem, several obstacles appear.
First, the time-inhomogeneous analogue of the semigroup, a so-called forward evolution system, no longer commutes with the (time-dependent) generator.
Second, the duality functions are then also naturally time dependent.
Therefore, the connection between ``generator duality'' and ``process duality'' becomes more delicate, because the time derivative of the duality function appears in the infinitesimal version of duality.
In particular, one sees that the process running forward in time has a dual running backward in time.
Therefore, the semigroup duality has to be replaced by duality between a forward and backward evolution system.
Finally, in the time-homogeneous setting, in the presence of a reversible (detailed balance) measure, one has an associated so-called cheap duality which follows directly from the detailed balance relation (\cite{redig2026duality} Chapter 2).
In the time-inhomogeneous setting it is no longer natural to have a stationary measure.
Rather, one has a trajectory of probability measures following the Kolmogorov forward equation, which we call a consistent trajectory of probability measures.
From this consistent trajectory of probability measures one can make a backward evolution system, using Bayes' rule.
The connection between the forward and backward system then replaces the notion of cheap duality coming from detailed balance.

The aim of this paper is thus to make first steps in developing a theory of duality for time-inhomogeneous Markov processes, and to illustrate it with examples from interacting particle systems and diffusion processes.
We are hereby inspired by the three particle systems which form the basis of \cite{redig2026duality}, namely the SEP, symmetric inclusion process (SIP) and independent random walkers.
For independent walkers, the extension is still quite complete, because there is a consistent trajectory of Poissonian product measures.
We also prove that the system of independent walkers has time-dependent symmetries which are weighted time-dependent sums of creation and annihilation operators.

In the interacting example (SIP, SEP), we prove that there are no longer proper time-dependent weighted sums of creation and annihilation operators which produce time-dependent symmetries. 
In order to obtain a consistent trajectory of (Gamma) product measures for the SIP, we add well-chosen time-dependent reservoirs.
In that setting we can then also prove a form of single particle duality.

The rest of our paper is organized as follows.
Section~\ref{sec:for_back_evol_sys} extends the notions of semigroups and generators to the time-inhomogeneous setting and provides some examples of those.
It also provides the definition of the replacement for detailed balance measures in the time-inhomogeneous setting.
In Section~\ref{sec:evol_sys_and_gen_duality}, we provide the definition of evolution system and generator duality for time-inhomogeneous processes and provide simple examples.
The creation of dualities from cheap dualities in combination with symmetries and intertwiners is laid out in Section~\ref{sec:symmetries_intertwining}.
This section ends with examples in which the (non-)existence of time-inhomogeneous symmetries is explored for independent random walkers and the symmetric inclusion and exclusion processes.
The last section, Section~\ref{sec:examples}, provides more involved examples of time-dependent dualities.

\section{Forward and backward evolution systems} \label{sec:for_back_evol_sys}
In this section we define the extension of a semigroup and generator to the time-inhomogeneous setting.
We also give the relationship between them.

In the time-homogeneous setting, the Markov semigroup associated to a Markov process $\{X(t): \, t\geq 0\}$ on the state space $\Omega$ is given by
\[S(t)f(x) = \E[f(X(t)) \mid X_0=x], \ x\in\Omega,\]
for a class of measurable functions $f:\Omega\to\R$.
By the Markov property combined with the time-homogeneity, we have $\E[f(X(t)) \mid X(s)=x] = S(t-s)f(x)$, for $0 \leq s\leq t$.
The analogue of the Markov semigroup associated to a time-inhomogeneous Markov process $\{Y(t): \, t\geq 0\}$ on the state space $\Omega$ is a two-parameter family $\{U(s, t): \, 0 \leq s \leq t\}$, which we call an evolution system.
This evolution system has the probabilistic interpretation for a class of measurable functions
\[
	U(s, t)f(y) = \E[f(Y(t)) \mid Y(s) = y], \ y\in\Omega.
\]
We denote  by $\caC(\Omega)$ the class of functions on which
these operators act.
Typically these are continuous functions, bounded continuous functions, or continuous functions vanishing at infinity.
Below we give the precise definition of such an evolution system \cite{bottcher2014feller}.
\begin{definition}[Forward evolution system] \label{def:evol_system}
	A family of linear operators $\{U(s, t): \, 0 \leq s \leq t \}$ on $\caC(\Omega)$ is a forward evolution system if it satisfies the following properties.
	\begin{enumerate}
		\item Identity: $U(s, s) = I$ for all $s \geq 0$, where $I$ denotes the identity operator.
		\item Normalization: $U(s, t) 1=1$.
		\item Positivity: $U(s, t) f \geq 0$ for $f \geq 0$ and $f \in \caC(\Omega)$.
		\item Contraction: $\|U(s, t) f\|_{\infty} \leq \|f\|_{\infty}$ for $f \in \caC(\Omega)$.
		\item Evolution property: $U(s, t)= U(s, r) U(r, t)$ for $s \leq r \leq t$.
		\item Strong continuity: $\lim_{\substack{(s, t) \to (u, w) \\ 0 \leq s \leq t \\ 0 \leq u \leq w}} \|U(s, t) f - U(u, w)f\|_{\infty} = 0$ for $f \in \caC(\Omega)$.
	\end{enumerate}
\end{definition}
In the time-homogeneous setting, the semigroup is the exponential of the generator (where the exponential has to be interpreted in the appropriate sense, i.e., via the Hille-Yosida theorem).
In the time-inhomogeneous setting this no longer holds.
The generator has to be replaced by a family of generators indexed by time, i.e., a time-dependent generator.
The relation between the evolution system and the time-dependent generator is given through a time derivative.
Additionally, the connection between the generator and the evolution system becomes more complex compared with the time-homogeneous setting, because the evolution system and the time-dependent generator do not commute in general.
Below we provide the definition of the time-dependent generator.
\begin{definition}[Time-dependent generator]
	Let $\{U(s,t): \, 0 \leq s \leq t\}$ be a forward evolution system.
	The time-dependent generator $L_s$ is defined by
	\begin{equation}\label{eq:def_tim_dep_gen}
		L_s f := \lim_{h \downarrow 0} \frac{U(s, s+h)f - f}{h},
	\end{equation}
	for all $f$ in the domain $\mathcal{D}(L_s)$ for which this limit exists in $\|\cdot\|_\infty$.
	The generator $L_s$ is related to the evolution system via the evolution equations
	\begin{align}\label{forwardevolutionequation}
		\begin{split}
			\pdert U(s,t) &= U(s,t) L_t, \\
			\pders U(s,t) &= -L_s U(s,t).
		\end{split}
	\end{align}
\end{definition}

\begin{remark}
	\sloppy Equation~\eqref{eq:def_tim_dep_gen} defines the generator through the right derivative (this is called the ``right generator'').
	Analogously, one can define the ``left generator'' by $L_{s}^{-} f:=\lim _{h \downarrow 0} \frac{U(s-h, s) f-f}{h}$.
	Under continuity conditions, see \cite{bottcher2014feller}, the left and right generators are identical.
	In our setting, all results hold for both the left and right generator, so without loss of generality we consider only the right generator.
\end{remark}
In our study of duality in the time-inhomogeneous setting, we naturally encounter processes evolving backward in time.
More precisely, in this setting the original process evolves forward in time, while the dual process generically evolves backward in time.
Therefore, we introduce a \textit{backward evolution system}, $\{\hat{U}(t, s): \, 0 \leq s \leq t \}$.
This family of operators satisfies all the defining properties of a forward evolution system, except that the evolution property is replaced by the backward evolution property $\hat{U}(t, s)= \hat{U}(t, r) \hat{U}(r, s)$ for $0 \leq s \leq r \leq t$.
A typical way to obtain a backward evolution system is via the backward flow of a process, i.e., a process $\{\hat{Y}(r): \, t\geq r\geq s\}$ evolving backward in time from time $t$ to time $s\leq t$.
So $\hat{U}(t, s)f(y) = \E[f(\hat{Y}(s)) \mid \hat{Y}(t) = y]$, with $0 \leq s\leq t$.
The time-dependent generator $\hat{L}_t$ of a backward evolution system is defined similarly to $L_t$.
However, the evolution equations differ from \eqref{forwardevolutionequation}, namely,
\begin{align*}
	\begin{split}
		\pdert \hat{U}(t, s) & = \hat{L}_{t} \hat{U}(t, s), \\
		\pders \hat{U}(t, s) & =- \hat{U}(t, s) \hat{L}_{s}.
	\end{split}
\end{align*}

Note that in the time-homogeneous setting, under standard conditions there is a one-to-one relation between the generator and the semigroup.
In the time-inhomogeneous setting, one generator can be related to two evolution systems, namely a forward and a backward evolution system.
This is possible because a generator and an evolution system do not necessarily commute.
The relation between the forward and backward systems when a generator and an evolution system commute is given in the following proposition.

\begin{proposition} \label{prop:commutation_gen_and_evolsys}
    Let $\{U(s, t): \, 0 \leq s \leq t\}$ be a forward evolution system.
	Let $\{L_t: \, t\geq 0\}$ be the corresponding time-dependent generator.
	\begin{itemize}
		\item [a)] Define
		      \begin{equation}\label{eq:backward_es_is_forward_es}
			      \hat{U}(t, s) = U(s, t).
		      \end{equation}
		      Assume that $\{\hat{U}(t, s): \, 0 \leq s \leq t\}$ is a backward evolution system.
		      Then we have the commutation of the forward evolution system with the time-dependent generator, i.e., for all $0 \leq s \leq t$,
		      \begin{equation} \label{eq:commutation_es_gen}
			      L_t U(s, t) = U(s, t) L_t.
		      \end{equation}
		\item[b)] If $\hat{V}(t,s)$ is a backward evolution system with time-dependent generator $L_t$, i.e., if  for all $0\leq s\leq t$ we have
		      \begin{equation}\label{veqo}
			      \pdert \hat{V}(t, s)  = L_t \hat{V}(t, s),
		      \end{equation}
		      and if \eqref{eq:commutation_es_gen} holds, then
		      $\hat{V}(t, s)= U(s,t)$.
	\end{itemize}
\end{proposition}
\begin{proof}
	We start with the proof of item a).
	Fix $0\leq s\leq r\leq t$.
	From the fact that \eqref{eq:backward_es_is_forward_es} defines a backward evolution system, we derive
	\[
		\hat{U}(t,r) \hat{U}(r, s)= U(r,t) U(s, r)=\hat{U}(t,s)= U(s,t).
	\]
	Using that $U(s,t)$ is a forward evolution system, we obtain
	\begin{equation}\label{bonkis}
		U(r,t) U(s, r)= U(s, r) U(r,t).
	\end{equation}
	To conclude \eqref{eq:commutation_es_gen}, take the derivative w.r.t. $t$ in \eqref{bonkis} and put $r=t$.

	To prove item b), assume \eqref{eq:commutation_es_gen}.
	Then we find
	\begin{equation*}
		\pdert U(s, t) = L_t U(s, t) = U(s, t) L_t.
	\end{equation*}

	Combined with \eqref{veqo}, we conclude that $U(s, t)$ and $\hat{V}(t, s)$ satisfy the same differential equation.
	Also the initial condition is the same because $U(s, s) = \hat{V}(s, s) = I$, therefore $U(s, t) = \hat{V}(t, s)$.
\end{proof}
An evolution system induces an evolution on probability measures (which with slight abuse of notation we still denote by $U(s,t)$) via
\[
\int U(s,t) f(x) \diff \mu_s (x)= \int f(x) \diff [ \mu_s U(s,t)](x).
\]
Probabilistically, $\mu_s U(s,t)$ is the distribution of $X(t)$ when started at time $s$ from
distribution $\mu_s$.
Below we make the evolution of the probability measures precise.
\begin{definition}[Consistent trajectory of probability measures]
	Let $\{U(s,t): \, 0 \leq s \leq t\}$ be a forward evolution system.
	A family of probability measures $\{\mu_s: \, s \geq 0\}$ is a consistent trajectory of probability measures for this evolution system if $\mu_s U(s,t) = \mu_t$ for all $0 \leq s \leq t$.
\end{definition}
A simple way to produce a consistent trajectory of probability measures is to take the distribution of $X(t)$ when initially started from $\mu$, i.e., $\mu_t= \mu U(0,t)$.
In the following subsection we provide more examples of consistent trajectories of probability measures.
\subsection{Examples}\label{subsec:examples_time_inhomo_gen}
Below we first give some examples of time-inhomogeneous generators and their evolution systems.
In the second half, we also provide a few examples of a consistent trajectory of probability measures.
\subsubsection*{Generators and the corresponding evolution system}
\subsubsection*{Example 2.1: Time-dependent Ornstein-Uhlenbeck process}
	      We consider the Ornstein-Uhlenbeck process which solves
	      \be\label{OUD}
	      \diff X(t) = - a(t) X(t) \diff t + \diff W(t),
	      \ee
	      where $a: [0,\infty)\to [0,\infty)$ is
	      a continuous function and $W(t)$ is Brownian motion.
	      Then the time-inhomogeneous generator corresponding to this process is
	      \[L_t = -a(t) x \frac{\diff}{\diff x} + \frac{1}{2}\frac{\diff^2}{\diff x^2},\]
	      and the evolution system is
	      \[
		      U(s, t) f(x) = \E\left(f\left( x e^{- \int_s^ta(u) \diff u} + \sigma(s, t) Z \right)\right),
	      \]
	      where $\sigma^2(s, t) = \int_s^t e^{-2 \int_r^t a(u) \diff u} \diff r. $
\subsubsection*{Example 2.2: Generator multiplied with time-dependent factor} \label{gen_mult_time_dep_factor}
	      We consider a time-inhomogeneous generator, which can be written as
	      \[
		      L_t = c(t) L,
	      \]
	      where $c(t)$ is some function of time and $L$ is a time-homogeneous generator.
	      This corresponds to a deterministic time change.
	      The evolution system is then given by
	      \[
		      U(s, t) = e^{\left(\int_s^t c(r) \diff r \right) L}.
	      \]
	      Notice that in this case, the evolution system commutes with the time-dependent generator.

Below we provide examples of consistent trajectory of measures for random walks and the Ornstein-Uhlenbeck process.
For the random walk consider a jump process on a countable set with the following generator
\be \label{eq:def_rw_gen}
L_t f(x) = \sum_{y } r_t(x, y) [f(y)-f(x)],
\ee
where $r_t(x, y)$ are time-dependent edge rates.
Then the adjoint of this generator (w.r.t. counting measure)  reads
\[
	L_t^{\star} \mu(x) = \sum_{y } [\mu(y) r_t(y, x) - \mu(x) r_t(x, y) ].
\]
A consistent trajectory of probability measures is then a solution of the Kolmogorov forward equation which is given by
\[
	\pdert \mu_t(x) = L_t^{\star} \mu_t(x).
\]
\subsubsection*{Consistent families of probability measures}
\subsubsection*{Example 2.3: Random walk on two sites}
	      Consider a random walk on the set $\{0, 1\}$, with generator given in \eqref{eq:def_rw_gen}.
	      We choose the following transition rates
	      \begin{align*}
		      r_t(0, 1) &= \nu(t), \\
		      r_t(1, 0) &= \lambda(t).
	      \end{align*}
	      Since $\mu_t$ is a probability measure, we have $\mu_t(0) + \mu_t(1) = 1$ for all $t \geq 0$.
	      Using this, we calculate the derivative of $\mu_t(0)$
	      \begin{align*}
		      \pdert \mu_t(0) &= r_t(1, 0) \mu_t(1) - r_t(0, 1) \mu_t(0) \\
		      &= \lambda(t) \mu_t(1) - \nu(t) \mu_t(0) \\
		      &= \lambda(t) (1 - \mu_t(0)) - \nu(t) \mu_t(0) \\
		      &= \lambda(t) - \left(\lambda(t) + \nu(t)\right) \mu_t(0).
	      \end{align*}
	      Then
	      \[
		      \mu_t(0) = \left[\mu_0(0) + \int_0^t \lambda(s) e^{\int_0^s \lambda(r) + \nu(r) \diff r} \diff s \right] \exp \left(-\int_0^t (\lambda(s) + \nu(s)) \diff s \right)
	      \]
	      is a solution to this differential equation and therefore defines a consistent trajectory of probability measures.
	\subsubsection*{Example 2.4: Symmetric time-inhomogeneous random walk on $\Z$}
	      Let us consider the random walk with time-dependent transition rates given by
	      \[
		      r_t(x, x+1) = r_t(x+1, x)  = r_t(x, x-1) = r_t(x-1, x) = \lambda(t).
	      \]
	      Then the solution of the following differential equation
	      \[\pdert \mu_t(x) = \lambda(t) \left(\mu_t(x+1) + \mu_t(x-1) - 2 \mu_t(x) \right)\]
	      provides a consistent trajectory of probability measures.
	      The solution is given by
	      \[
		      \mu_t(x)= \E_x^{\RW} \left(\mu_0(X(\tau(t)))\right),
	      \]
	      where $\{X(t): \, t\geq 0 \}$ is the simple continuous-time random walk with rates $r(x,x+1)= r(x,x-1)=1$, $\E_x^{\RW}$ denotes expectation w.r.t. this random walk starting at $X(0)=x$ and $\tau(t)= \int_0^t \lambda(s) \diff s$.
	\subsubsection*{Example 2.5: Ornstein-Uhlenbeck process} 
	      For the Ornstein-Uhlenbeck process given in \eqref{OUD}, the time-dependent normal distribution
	      \[
		      \mu_t= \caN\left(0, \si^2 e^{-2K(0,t)}+ \int_0^t e^{-2K(s,t)}\diff s\right),
	      \]
	      where $K(s,t)= \int_s^t a(r) \diff r$ is an example of a consistent trajectory of probability measures.

\section{Evolution system duality and generator duality}\label{sec:evol_sys_and_gen_duality}
In this section we define duality for time-inhomogeneous Markov processes and derive some basic properties.
It is natural to consider duality functions which besides depending on the state space variables of the two processes (the process and its dual) also depend explicitly on time.
We define the analogue of semigroup duality and generator duality and the relation between them.
In the time-inhomogeneous setting, it becomes clear that a meaningful definition of duality is necessarily between a forward and a backward evolution system.
In the time-homogeneous setting the backward evolution system is not ``visible'' because the semigroup and the generator commute.

In this section, if we want to emphasize the connection between $U$ and the corresponding process $X=\{X(t): \, t\geq 0\}$, we use the notation $U^X(s,t)$.
We denote by $X(s,x;t)$ the process $X$ at time $t$ when started at time $0 \leq s\leq t$ from $x$, in particular $X(s, x; s) = x$.
The relation with the forward evolution system then reads $U^X(s, t) f(x) = \E(f(X(s, x;t)))$.
Note that this expectation is a conditional expectation.
For a backward process, $\{\hat{Y}(t): \, t \geq 0 \}$, associated to a backward evolution system, we use the following notation for the evolution system and generator: $\hat{U}^Y(t, s)$ and $\hat{L}_t^Y$.
Correspondingly, we write $\hat{Y}(t, y; s)$ for the state at time $s$ of a backward process $\hat{Y}$ starting from $y$ at time $t\geq 0$.
These random variables $\hat{Y}(t, y; s)$ are related to the backward evolution system via $\hat{U}^Y(t, s) f(y)= \E(f(\hat{Y}(t, y; s)))$.
Below we define duality for evolution systems.

\begin{definition}[Evolution system duality]
	\sloppy Let $\{U^X(s, t): \, 0 \leq s \leq t\}$ and $\{\hat{U}^Y(t, s): \, 0 \leq s \leq t\}$ denote a forward and backward evolution system, respectively.
    Let $\Omega$ and $\hat{\Omega}$ be the state spaces of the processes with evolution systems $U^X$ and $\hat{U}^Y$, respectively. 
	Let $D: [0,\infty)\times \hat{\Omega} \times \Omega \to \R$ denote a measurable function.
	We say that $D$ is a duality function
	between $U^X(s, t)$ and $\hat{U}^Y(t, s) $ if for all $x \in \Omega, \, y \in \hat{\Omega}$,  and for all $0\leq s\leq t$ we have
	\begin{equation} \label{eq:evol_sys_duality}
		U^X(s, t)D(t, y, \cdot)(x) = \hat{U}^Y(t, s) D(s, \cdot, x)(y).
	\end{equation}
\end{definition}
\br
In terms of expectations, \eqref{eq:evol_sys_duality}
reads
\[
\E( D(t, y, X(s,x;t)))= \E( D(s, \hat{Y}(t,y;s), x)).
\]
Therefore, we can think of this relation as duality between the forward process $\{X(r): \, 0 \leq s\leq r\leq t\}$, and
the backward process $\{ \hat{Y}(r): \, 0 \leq s\leq r\leq t\}$.
\er

In the time-homogeneous setting, there is a symmetry between the original process and the dual process in the sense that it does not make a difference which of the two processes is called the original or dual process.
In the time-inhomogeneous setting, this symmetry is no longer present.
This is because the duality function is time-dependent, and because $X$ has a forward and $\hat{Y}$ a backward evolution system.

If the duality functions are in the domain of the infinitesimal generator and are differentiable, then generator duality can be defined.
\begin{definition}[Time-dependent generator duality]
	Let $\{L_t^X: \, t\geq 0\}$ denote the generator of the process $\{X(s, \cdot \,; t): \, 0 \leq s \leq t\} $ and $\{\hat{L}_t^Y: \, t \geq 0\}$ the generator of the backward process $\{\hat{Y}(t, \cdot \, ; s): \, 0 \leq s \leq t\}$.
	Let $\Omega$ and $\hat{\Omega}$ be the state spaces of the processes $X$ and $\hat{Y}$, respectively. 
    Let $D: [0, \infty) \times \hat{\Omega} \times \Omega \to \R$ be a function such that $D(t, y, \cdot)$ is in the domain of $L_t^X$ and $D(t, \cdot, x)$ is in the domain of $\hat{L}_t^Y$, and $D$ is differentiable in $t \geq 0$.   
	We then say that $D$ is a duality function between $L_t^X$ and $\hat{L}_t^Y$ if for all $x \in \Omega, \, y \in \hat{\Omega}$, we have
	\begin{equation*}
		L_t^X D(t, y, \cdot)(x) + \pdert D(t, y, x) = \hat{L}_t^Y D(t, \cdot, x)(y).
	\end{equation*}
\end{definition}

\begin{remark}
	In the following, we drop the notation involving the dot.
	Instead, the superscript indicates on which variable the evolution system or generator acts.
\end{remark}

The following theorem extends the equivalence relation between generator and semigroup duality of \cite{redig2026duality} to time-inhomogeneous processes where the duality function depends on time.

\begin{theorem}[Relation between generator and evolution system duality] \label{thm:relation_gen_dual_evol_dual}
	Let $\{U^X(s, t): \, 0 \leq s \leq t\}$ be an evolution system and let $\{\hat{U}^Y(t, s): \, 0 \leq s \leq t\}$ denote a backward evolution system of the processes $\{X(s, \cdot \,; t): \, 0 \leq s \leq t\}$ and the dual process $\{\hat{Y}(t, \cdot \, ; s): \, 0 \leq s \leq t\}$, respectively.
	Let $\{L_t^X: \, t\geq 0\}$ and $\{\hat{L}_t^Y: \, t\geq 0\}$ denote the generators corresponding to the forward and backward evolution system.
    Let $\Omega$ and $\hat{\Omega}$ be the state spaces of the processes $X$ and $\hat{Y}$, respectively.
	We assume that $D: [0, \infty) \times \hat{\Omega} \times \Omega \rightarrow \R$ is such that for all $x \in \Omega$ and for all $t \geq 0$, $D(t, \cdot, x)$ is in the domain of $\hat{L}_t^Y$.
    Additionally, for all $y \in \hat{\Omega}$ and for all $t \geq 0$, $D(t, y,\cdot)$ is in the domain of $L_t^X$. 
    Further assume that $D$ is in $C^1$ for $t \geq 0$.
	\begin{itemize}
		\item[a)] Evolution system duality implies generator duality, i.e., if for all $x \in \Omega, \, y \in \hat{\Omega}$ and $0 \leq s \leq t$,
		      \begin{equation}\label{eq:evol_sys_duality_time_dep_dual_func}
			      U^X(s, t)D(t, y, x) = \hat{U}^Y(t, s)D(s, y, x),
		      \end{equation}
		      then we have for all $t \geq 0$
		      \begin{equation}\label{eq:gen_duality_time_dep_dual_func}
			      L_t^X D(t, y, x) + \frac{\partial}{\partial t} D(t, y, x) = \hat{L}_t^Y D(t, y, x).
		      \end{equation}
		\item[b)] In the other direction, if for all $0 \leq s \leq t$, $\hat{U}^Y(t, s) D(t, y, \cdot)$ is in the domain of $L_t^X$, then generator duality, i.e., \eqref{eq:gen_duality_time_dep_dual_func} implies evolution system duality, i.e., \eqref{eq:evol_sys_duality_time_dep_dual_func}.
	\end{itemize}
\end{theorem}
\begin{proof}
	We first prove item a).
	Assume \eqref{eq:evol_sys_duality_time_dep_dual_func}.
	Then
	\begin{equation*}
		\pdert U^X(s, t) D(t, y, x) = U^X(s, t) L_t^X D(t, y, x) + U^X(s, t) \pdert D(t, y, x).
	\end{equation*}
	We put $s=t$ to obtain
	\begin{equation*}
		\pdert U^X(s, t) D(t, y, x) \bigg|_{s=t} = L_t^X D(t, y, x) + \pdert D(t, y, x).
	\end{equation*}
	Furthermore,
	\begin{equation*}
		\pdert \hat{U}^Y(t, s)D(s, y, x) = \hat{L}_t^Y \hat{U}^Y(t, s)D(s, y, x).
	\end{equation*}
	Putting $s=t$ gives
	\begin{equation*}
		\pdert \hat{U}^Y(t, s) D(s, y, x) \bigg|_{s=t} = \hat{L}_t^Y D(t, y, x).
	\end{equation*}
	Second, we prove item b).
	We consider the interpolating quantity
	\begin{equation*}
		\psi(r) := U^X(s, s+r) \hat{U}^Y(t, s+r) D\left(s+r, y, x\right).
	\end{equation*}
	Then we have that $\psi(0)$ equals the right hand side of \eqref{eq:evol_sys_duality_time_dep_dual_func} and $\psi(t-s)$ equals the left hand side of \eqref{eq:evol_sys_duality_time_dep_dual_func}.
	Therefore we want to derive $\psi(0)=\psi(t-s)$ from \eqref{eq:gen_duality_time_dep_dual_func}.
	Taking the derivative w.r.t. $r$, assuming \eqref{eq:gen_duality_time_dep_dual_func}, and using that a generator and evolution system acting on different variables ($x$, respectively $y$) commute, we obtain
	\begin{equation*}
		\frac{\diff}{\diff r} \psi(r) = U^X(s, s+r) \hat{U}^Y(t, s+r) \left( \left[L_{s+r}^X - \hat{L}_{s+r}^Y + \frac{\partial}{\partial r} \right] D\left(s+r, y, x\right) \right)
		=0.
	\end{equation*}
	The above calculation implies $\psi(0)=\psi(t-s)$.
\end{proof}
\begin{remark}
	Note that we did not use the normalization condition in the proof of Theorem~\ref{thm:relation_gen_dual_evol_dual}.
	Therefore, this theorem also holds for evolution systems that are not necessarily mass preserving, for instance processes with a Feynman-Kac factor.
\end{remark}
The following corollary states what happens when the forward and backward processes are the same.
\begin{corollary}\label{cor:eqv_commutation_and_forward_equals_backward}
	\sloppy Let $\{U^X(s, t): \, 0 \leq s \leq t\}$ and $\{U^Y(s, t): \, 0 \leq s \leq t\}$ be two evolution systems of the processes $\{X(s, \cdot; t): \, 0 \leq s \leq t \}$ and $\{Y(s, \cdot; t): \, 0 \leq s \leq t \}$, respectively. 
    Let $\{L_t^X: \, t \geq 0\}$ and $\{L_t^Y: \, t \geq 0\}$ denote the generators corresponding to the evolution processes. 
    Let $\Omega$ and $\hat{\Omega}$ be the state space of the processes $X$ and $Y$, respectively.
    We assume that $D: [0, \infty) \times \hat{\Omega} \times \Omega \rightarrow \R$ is such that for all $x \in \Omega$ and for all $t \geq 0$, $D(t, \cdot, x)$ is in the domain of $L_t^Y$.
    Additionally, for all $y \in \hat{\Omega}$ and for all $t \geq 0$, $D(t, y,\cdot)$ is in the domain of $L_t^X$ and for all $y \in \hat{\Omega}$ and for all $0 \leq s \leq t$, $U^Y(s, t) D(s, y, \cdot)$ is in the domain of $L_t^X$. 
    Further assume that $D$ is differentiable in $t \geq 0$.
    Then for all $x \in \Omega$ and for all $y \in \hat{\Omega}$ and for all $ 0 \leq s \leq t$
	\be \label{eq:evol_duality}
		U^X(s, t) D(t, y, x) = U^Y(s, t) D(s, y, x)
	\ee
    implies that for all $t \geq 0$
    \be \label{eq:gen_duality}
		L_t^X D(t, y, x) + \pdert D(t, y, x) = L_t^Y D(t, y, x).
	\ee
	If in addition the generator and the evolution system of the process $Y$ commute for all $0 \leq s \leq u \leq v \leq t$, i.e., 
	\begin{equation} \label{eq:commutation}
		L^Y_v U^Y(u, v) = U^Y(u, v) L^Y_v,
	\end{equation}
    then \eqref{eq:gen_duality} implies \eqref{eq:evol_duality}.
    In particular, if $D$ does not depend on time, then for all $0 \leq s \leq t$ 
	\be \label{eq:evol_duality_2}
		U^X(s, t) D(y, x) = U^Y(s, t) D(y, x)
	\ee
    implies that for all $t \geq 0$
	\be \label{eq:gen_duality_2}
		L_t^X D(y, x) = L_t^Y D(y, x).
	\ee
	If in addition \eqref{eq:commutation} holds, then \eqref{eq:gen_duality_2} implies \eqref{eq:evol_duality_2}.
\end{corollary}
\begin{proof}
	Notice that $U^Y(s, t)$ is a backward evolution system with generator $L_t^Y$. 
    The proof of this equivalence is similar to that of Theorem~\ref{thm:relation_gen_dual_evol_dual}.
\end{proof}

\subsection{A few simple examples}
Here we provide examples of time-inhomogeneous processes and their duality functions.
In two of these examples the backward process is the same as the forward process.

\subsubsection*{Example 3.1: Time-inhomogeneous Wright-Fisher diffusion}
We consider the diffusion process
\begin{equation*}
	\diff X(t) = \sigma_t \sqrt{X(t)(1-X(t))} \diff W(t),
\end{equation*}
where $W(t)$ denotes Brownian motion. 
The  time-dependent generator is given by
\begin{equation*}
	L^X_t = \frac{1}{2} \sigma^2_t x(1-x) \frac{\partial^2}{\partial x^2}.
\end{equation*}
We have $L^X_t = \sigma^2_t L^X$, with $L^X$ the generator of the time-homogeneous Wright-Fisher diffusion, namely
\begin{equation*}
	L^X f(x) = \frac{1}{2} x(1-x) \frac{\partial^2}{\partial x^2} f(x).
\end{equation*}
As a consequence, the evolution system
can be written as follows
\begin{equation*}
	U^X(s, t) f(x) = e^{T(s,t) L^X} f(x)=
	\E\left( f\left(Z(T(s,t))\right)\right),
\end{equation*}
where $T(s,t)= \int_s^t \sigma_r^2 \diff r $ and $\{Z(t): \, t\geq 0\}$ denotes the time-homogeneous Wright-Fisher diffusion with generator $L^X$.

From the time-homogeneous setting it is known that the Wright-Fisher process is dual to the Kingman's coalescent block counting process.
The time-homogeneous generator of the latter process
$\{N(t): \, t\geq 0\}$ is given by
\begin{equation*}
	L^N f(n) = \frac{1}{2} n(n-1) (f(n) - f(n-1)).
\end{equation*}
Using the duality function $D(n, x) = x^n$, we have $L^XD(n, x) = L^N D(n, x)$.
The same duality function gives us the (time-dependent) generator duality
\begin{equation*}
	L_t^X D(n, x) = \sigma_t^2 L^X D(n, x) = \sigma_t^2 L^N D(n, x) = L_t^N D(n, x),
\end{equation*}
with $L_t^N f(n) = \sigma_t^2 \frac{1}{2} n(n-1) (f(n) - f(n-1))$.
As a consequence, we have evolution system duality which then reads
\[
	U^X(s,t) D(n,x)= \hat{U}^N(t,s) D(n, x),
\]
where
\[
	\hat{U}^N(t,s) f(n)= \E_n \left(f\left( N(T(s,t))\right)\right)= U^N (s,t) f(n).
\]
Note that the dual process does not have to be reversed, because the time dependence in this example corresponds to a deterministic time change.
Therefore, the time-dependent generator and the evolution system commute.

\subsubsection*{Example 3.2:
	Ornstein-Uhlenbeck process}
Consider the Ornstein-Uhlenbeck process solving the SDE
\[
\diff X(t) = -a(t) X(t) \diff t+ \si_t \diff W(t),
\]
where $W(t)$ denotes Brownian motion. 
Then the time-dependent generator reads
\[
	L^X_t= -a(t) x\partial_x + \frac12 \si_t^2\partial_x^2.
\]
Consider
\[
	D(t,y,x) = e^{c(t) xy},
\]
then acting with $\partial_t + L^X_t$ on the $x$ and $t$ variables gives
\[
	(\partial_t + L^X_t)D(t,y,x)
	= \left(\left(\frac1{c(t)}\pdert c(t) - a(t)\right)y\partial_y  + \frac12 \si_t^2 c(t)^2 y^2 \right) D(t,y,x).
\]
Therefore, if
\be\label{difi}
	\left(\frac1{c(t)}\pdert c(t) - a(t)\right) =0,
\ee
then we find that the dual time-dependent generator is the multiplication operator
\[
	\hat{L}^Y_t f(y)= \frac12 \si_t^2 c(t)^2 y^2 f(y),
\]
which corresponds to the reversed multiplicative evolution system
\[
	\hat{U}^Y(t,s)f (y)= \exp\left(\frac{y^2}{2} \int_s^t \si^2_r c(r)^2 \diff r\right) f(y).
\]
Thus we obtain the duality
\[
	U^X(s,t) D(t,y, x)= \hat{U}^Y(t,s) D(s, y, x),
\]
which reads in probabilistic notation
\be\label{probnot}
	U^X(s, t) D(t, y, x) = \E\left( e^{ c(t) X(t) y}|X(s)= x\right)
	=
	e^{\frac{y^2}{2}\int_s^t \si^2_r c(r)^2\diff r} e^{c(s) xy }.
\ee
There is a  unique solution $c(r)$ with $s\leq r\leq t$ of
the differential Equation~\eqref{difi} with end condition $c(t)=1$.
Then, using \eqref{probnot} we conclude that, conditional on $X(s)=x$,
$X(t)$ is normally distributed  with expectation
$c(s) x$ and variance
$\int_s^t \si^2_r c(r)^2\diff r$.

\subsubsection*{Example 3.3: Multiplication evolution system}
Define
\begin{equation*}
	\diff X(t) = \sigma_t \diff W(t),
\end{equation*}
where $W(t)$ denotes Brownian motion. 
We define the following two evolution systems
\begin{equation*}
	U^X(s, t) f(x)= \E(f(X(t)) \mid X(s) = x),
\end{equation*}
and
\begin{equation*}
	\hat{U}^Y(t, s) f(y)=\exp \left(y^2 \int_s^t \frac{1}{2} \sigma_r^2 \diff r\right) f(y).
\end{equation*}

Using the duality function $D(t, y, x) = D(y, x)=e^{y x}$, we find the evolution system duality
\begin{equation*}
	U^X(s, t) D(y, x) = \hat{U}^Y(t, s) D(y, x).
\end{equation*}
A direct computation gives
\begin{align*}
	U^X(s, t) e^{yx} &= \E \left[e^{y\left [x + \int_s^t \sigma_r \diff W(r) \right]} \right] \\ &= e^{yx} e^{\frac{y^2}{2}\int_s^t \sigma_r^2 \diff r} \\ &= \hat{U}^Y(t, s) e^{yx}.
\end{align*}

On the other hand, we can also derive the generator duality.
The generator of the evolution system $U^X$ is given by
\begin{equation*}
	L^X_t f(x) = \frac{1}{2} \sigma_t^2 \frac{\diff^2}{\diff x^2} f(x).
\end{equation*}
The generator of the evolution system $\hat{U}^Y$ is given by
\[
	\hat{L}^Y_t f(y) = \frac{1}{2} \sigma_t^2 y^2 f(y).
\]
It is then immediate that
\[
	L^X_t e^{yx} = \hat{L}^Y_t e^{yx}.
\]
Since $\hat{L}^Y_t$ and $\hat{U}^Y(t, s)$ commute, by Proposition~\ref{prop:commutation_gen_and_evolsys} it is clear that $\hat{U}^Y(t, s) = U^Y(s, t)$.

\section{Dualities, symmetry and intertwining}\label{sec:symmetries_intertwining}
In the time-homogeneous case, it is well-known that a reversible measure, equivalently, a detailed balance measure, automatically gives a so-called cheap duality function, see \cite{redig2026duality}, Section 2.3.
Acting with symmetries on this cheap duality function produces new duality functions.
In the time-inhomogeneous setting, it is no longer natural to work with a detailed balance measure.
Rather, we define a cheap duality via a consistent trajectory of probability measures, and the associated time-reversed process.
Also, the notions of symmetry and intertwining have to be modified to incorporate the behavior across time. 
\subsection{Notation and definitions}\label{subsec:from_c.m.t_to_duality}
To explain the concepts of cheap duality, symmetry and intertwining in the easiest context, we mostly restrict here to finite (or countable) state spaces.
So we consider a time-inhomogeneous Markov chain on the finite state space $S$.
We denote the transition probabilities of the Markov chain by $p_{s,t}(x,y)= \mathbb{P} (X(t)=y|X(s)=x)$.
The evolution systems $\{U(s,t): \, 0 \leq s \leq t\}$ are then given by
\[
	U(s,t) f(x)= \sum_y p_{s,t} (x,y) f(y).
\]
This family of operators satisfies the evolution property (cf. Definition~\ref{def:evol_system})
\be\label{semi}
U(s,r)U(r,t) f = U(s,t) f,
\ee
for $0 \leq s\leq r\leq t$.
The relation~\eqref{semi} corresponds to the
Chapman-Kolmogorov equation
\[
\sum_y p_{s,r}(x,y) p_{r,t} (y,z)= p_{s,t} (x, z)
\]
for all $x,z\in S$ and $0 \leq s\leq r\leq t$.
In this setting, the notion of a consistent trajectory of probability measures simplifies to a collection
$\{\mu_t: \, t\geq 0\}$ satisfying
$\mu_t(x)\geq 0, \sum_x \mu_t(x)=1$, and
\[
\sum_y \mu_s (y) p_{s,t}(y,x) =\mu_t(x),
\]
for all $0 \leq s\leq t, \, x\in S$.
The infinitesimal transition rates are defined via
\[
\lim_{\delta\to 0} \frac{p_{t,t+\delta}(x,y)-\delta_{x,y}}{\delta}= r_t(x,y).
\]
The time-dependent generator then reads
\[
L_t f(x)= \sum_y r_t(x,y)( f(y)-f(x)).
\]
So $L_t$ is a matrix where the rows sum to zero.
Its matrix elements are given by the transition rates $L_t(x, y)= r_t(x,y)$, and $L_t(x,x)= -\sum_{y\not= x} r_t(x,y)$.

We say that the rates $r_t(x,y)$ correspond to an irreducible chain if
for all $x,y\in S$ there exists a path $x_0=x, \ldots, x_n=y$ such
that $r_t(x_k,x_{k+1})>0$ for all $k\in \{0,\ldots, n-1\}$. In this section we will assume that for all $t\geq 0$ the rates $r_t(x,y)$ correspond to an irreducible chain. This implies that starting from an initial measure $\mu_0>0$,
for all $t>0$, it holds $\mu_t>0$.
A consistent trajectory of measures satisfies the Kolmogorov forward equation which reads
\[
\pdert \mu_t(x) = \sum_{y} (\mu_t(y)r_t(y,x)- \mu_t(x) r_t(x,y))= (\mu_t L_t^{\top})_x,
\]
where $\top$ denotes transposition, and $\mu_t$ is a row vector with entries $\mu_t(x), x\in S$.
With these notions, we define a reversed Markov process.
\begin{definition}[Reversed process]
	Let $\{\mu_t: \, t \in [0, \infty)\}$ be a consistent trajectory of probability measures with $\mu_t > 0$ for all $t \geq 0$. 
    Then we define a backward transition probability $\p_{t,s}$ for $0 \leq s \leq t$ via Bayes' formula
	\be\label{bayesbob}
		\p_{t,s}(x,y)= \frac{p_{s,t}(y,x) \mu_s(y)}{\mu_t(x)}.
	\ee
	To the backward transition probabilities corresponds the backward transition operator, i.e., for $0 \leq s \leq t$
	\[
		\U (t,s) f(x)= \sum_y \p_{t,s} (x,y) f(y)= \E(f(X^{\mathrm{rev}}(t,x;s))).
	\]
\end{definition}
\begin{remark}
	Since $\{\mu_t: \, t \geq 0 \}$ is a consistent trajectory of probability measures, $\p_{t,s}(x,y)$ is a proper backward transition probability, i.e.,
	\[
		\sum_y\p_{t,s}(x,y)=1,
	\]
	and for $0 \leq s \leq r \leq t$ we have the backward Chapman Kolmogorov equation:
	\[
		\sum_y\p_{t,r}(x,y)\p_{r, s}(y,z)= \p_{t,s} (x,z).
	\]
	As a consequence, the operators $\U (t,s)$ satisfy the backward evolution property.
\end{remark}
\br\label{choiceofmu}
\begin{itemize}
	\item[a)]The consistent trajectory of probability measures $\mu_t$ is not unique.
	      Therefore the backward transition probability depends on this choice of $\mu_t$.
        To alleviate notation, we do not explicitly write the dependence on the choice of $\mu_t$.
	      Note that this non-uniqueness even holds in the time-homogeneous case, i.e., we can have a time-inhomogeneous backward process corresponding to the time reversal of a homogeneous forward process.
    \item[b)] The assumption of irreducibility stated above ensures that the reversed process can be defined in the way described above, i.e., the quotient in \eqref{bayesbob} is well-defined.
	\item[c)] In case the transition rates are time-homogeneous (i.e., when $p_{s,t}(x,y)= p_{t-s}(x,y)$), and $\mu$ does not depend on time (i.e., is a stationary measure), then we have
	      \[
		      \p_{t,s}(x,y)=\tilde{p}_{t-s}(x,y)=
		      \frac{p_{t-s}(y,x)\mu(y)}{\mu(x)}.
	      \]
	      This is the transition probability of the backward process
	      (which is then also stationary and time-homogeneous), which  in the reversible (detailed balance) case coincides with $p_{t-s}(x,y)$.
    \end{itemize}
\er

We can think of $X(s,x; t)$ as the `stochastic flow' that flows forward in time from time $s$ to time $t$, and of $X^{\mathrm{rev}}(t,y; s)$ as the complementary backward stochastic flow that flows backward in time from time $t$ to time $s$.

In the following proposition we identify the reversed evolution system as the adjoint of the forward evolution system and also compute the generator of the backward evolution system.

\bp \label{prop:adjoint_U}
\ben
\item[a)] The generator of the backward evolution operator is given by
\begin{equation*}
	L_t^\rev f(x) = \sum_z \frac{\mu_t(z)}{\mu_t(x)} L_t(z, x) [f(z) - f(x)].
\end{equation*}
\item[b)] The backward evolution operator $\U(t,s)$
is the adjoint of the forward evolution operator $U(s,t)$ when the latter is viewed as mapping from
$L^2 (\mu_t)$ to $L^2 (\mu_s)$.
\item[c)] The adjoint of $L_t$ is given by $L_t^\rev + \frac{\pdert \mu_t}{\mu_t}$ when $L_t$ is viewed as a mapping from $L^2(\mu_t)$ to $L^2(\mu_t)$.
\een
\ep

\bpr
We first prove item a).
We start by noticing that
\[
	\pdert p_{s,t} (x,y)= \sum_z p_{s,t}(x,z) L_t(z,y).
\]
We use the above to compute $L_t^\rev$ by calculating the derivative of $U^\rev(t, s)$,
\begin{align*}
	\pdert U^\rev(t, s) f(x) &= \pdert \sum_y p_{t, s}^\rev (x, y) f(y) = \sum_y \pdert \frac{p_{s, t}(y, x) \mu_s(y)}{\mu_t(x)} f(y) \\
	&= \sum_y \frac{\left(\pdert p_{s, t}(y, x) \right) \mu_s(y) \mu_t(x) - p_{s, t}(y, x) \mu_s(y) \pdert \mu_t(x)}{\mu_t(x)^2} f(y) \\
	&= \sum_{y, z} \left( p_{s, t}(y, z) L_t(z, x) \frac{\mu_s(y)}{\mu_t(x)} - \frac{p_{s, t}(y, x) \mu_s(y) \mu_t(z) L_t(z, x)}{\mu_t(x)^2} \right) f(y) \\
	&= \sum_{y, z} \left(p^\rev_{t, s}(z, y) \frac{\mu_t(z)}{\mu_t(x)} L_t(z, x) - p^\rev_{t, s}(x, y) \frac{\mu_t(z)}{\mu_t(x)} L_t(z, x) \right) f(y) \\
	&= \sum_z \frac{\mu_t(z)}{\mu_t(x)} L_t(z, x)  [U^\rev(t, s)f(z)- U^\rev(t, s)f(x)] \\
	&= L^\rev_t U^\rev(t, s) f(x),
\end{align*}
with
\begin{equation*}
	L_t^\rev f(x) = \sum_z \frac{\mu_t(z)}{\mu_t(x)} L_t(z, x) [f(z) - f(x)].
\end{equation*}

To prove item b), we start by showing that $U(s, t)$ is a contraction between $L^2(\mu_t)$ and $L^2(\mu_s)$.
If $f \in L^2(\mu_t)$, then, using Jensen's inequality
\[
   \|U(s,t) f\|^2_{L^2(\mu_s)}= \int (U(s, t)f)^2 \diff \mu_s \leq \int U(s, t) (f^2) \diff \mu_s = \int f^2 \diff \mu_t= \|f\|_{L^2(\mu_t)}^2.
\]
which shows the claimed contraction property.
To see that $U^\rev(t, s)$ is the adjoint of $U(s, t)$, using \eqref{bayesbob} we obtain 
\begin{align*}
	\sum_{x,y} p_{s,t} (x,y) f(y) g(x) \mu_s(x)
	&=
	\sum_{y,x} f(y) \p_{t,s} (y,x) g(x) \mu_t (y),
\end{align*}
which can be rewritten as
\[
\int (U(s,t) f) g \diff \mu_s= \int (\U(t,s)g) f \diff \mu_t.
\]

To prove item c), notice that we can rewrite the inner product of $L_t f$ and $g$ as follows
\[
\sum_{x, y} L_t(x, y) f(y) g(x) \mu_t(x) = \sum_{x, y} f(x) \mu_t(x) \left[\frac{\mu_t(y)}{\mu_t(x)} L_t(y, x) g(y) \right].
\]
So the adjoint of $L_t$ is given by $L_t^\dagger g(x) = \sum_y \frac{\mu_t(y)}{\mu_t(x)} L_t(y, x) g(y)$.
Then we can write $L_t^\rev$ as follows
\begin{align*}
    L_t^\rev g(x) &= \sum_y \frac{\mu_t(y)}{\mu_t(x)} L_t(y, x) [g(y) - g(x)] \\
    &= L_t^\dagger g(x) - g(x) \frac{\sum_y \mu_t(y) L_t(y, x)}{\mu_t(x)} \\
    &= L_t^\dagger g(x) - g(x) \frac{\pdert \mu_t(x)}{\mu_t(x)}.
\end{align*}
Rewriting the above equation gives
\[
L_t^\dagger g(x) = L_t^\rev g(x) + g(x) \frac{\pdert \mu_t(x)}{\mu_t(x)}.
\]
Note that $L_t^\dagger$ is not a generator when $\mu_t$ truly depends on time. 
\epr


\subsection{Self-duality, symmetries and intertwining}
In this section we extend the time-homogeneous notion of a cheap duality and establish the notion of time-inhomogeneous symmetries and intertwiners.
\begin{definition}[Time-inhomogeneous self-duality]
	A time-inhomogeneous process $\{X(t): \, t \geq 0\}$ with evolution system $\{U(s, t): \, 0 \leq s \leq t\}$ and time-dependent generator $\{L_t: \, t \geq 0\}$ is self-dual to its reversed process $\{X^\rev(t): \, t\geq 0 \}$ with evolution system $\{U^\rev(t, s): \, 0 \leq s \leq t \}$ and time-dependent generator $\{L_t^\rev: \, t \geq 0\}$ with duality function $D: [0, \infty) \times S \times S \to \R$ if
	\[
		U(s, t) D(t, y, x) = U^{\rev}(t, s) D(s, y, x),
	\]
	where $U(s, t)$ acts on $x$ and $U^\rev(t, s)$ acts on $y$.
	Equivalently, if $D$ is differentiable in $t \geq 0$,
	\[
		L_t D(t, y, x) + \pdert D(t, y, x) = {L_t^\rev} D(t, y, x),
	\]
	where again $L_t$ acts on $x$ and $L^\rev_t$ acts on $y$.
\end{definition}

The following proposition gives self-duality between the process and its reverse using a cheap duality function built from a consistent trajectory of probability measures.
\bp \label{prop:cheap_duality}
Let $\{\mu_t: \, t\geq 0\}$ be a consistent trajectory of probability measures for the evolution system $\{U(s, t): \, 0 \leq s \leq t\}$.
Define for $t \geq 0$, $x, y \in S$ the cheap duality function
\[
	D_{\mathrm{cheap}}(t,y,x)= \frac{1}{\mu_t (x)} \delta_{x,y}.
\]
Then we have the following cheap duality between the process $\{X(t): \, t\geq 0\}$ and its reverse, $\{ X^\rev(t): \, t\geq 0\}$ with corresponding evolution system and generator for all $0 \leq s \leq t$ and $x, y \in S$
\[
	U(s,t) D_{\mathrm{cheap}}(t, y,x) = \U(t,s) D_{\mathrm{cheap}}(s, y, x),
\]
where $U(s, t)$ acts on $x$ and $U^\rev$ acts on $y$.
\ep
\bpr
We have
\begin{align*}
	\U(t,s) D(s, y, x) &=
	\sum_{y'}\p_{t,s}(y, y') D(s, y', x)\\
	&= \p_{t,s}(y,x)\frac{1}{\mu_s(x)}\\
	&= p_{s,t}(x,y) \frac{1}{\mu_t(y)}\\
	&= \sum_{x'}p_{s,t} (x,x') D(t, y,x')
\end{align*}
and the result follows.
\epr

The duality relation can be written in matrix form by defining $D(t)_{yx} = D(t, y, x)$, which gives
\[
D(t) U(s,t)^{\top} = \U(t,s) D(s),
\]
where $\top$ denotes transposition.
Also the infinitesimal version of the duality relation can be written in matrix form
\[
{L_t^{\rev}}
D(t)-D(t) L_t^{\top}=\dert D(t).
\]

\subsubsection{New dualities from symmetries}
In the time-homogeneous setting, non-trivial duality relations can be obtained by acting with symmetries of the process on duality functions.
To extend this, we first introduce symmetries in the time-inhomogeneous setting.
We introduce these notions for general backward processes, which are denoted by a hat.

\bd[Family of symmetries for evolution systems]
We call $\{T(t): \, t \geq 0\}$ a family of symmetries for
$\{U(s,t): \, 0 \leq s \leq t\}$ if
for all $0 \leq s\leq t$
\begin{equation}\label{eq:evol_sys_symmetry}
	U(s,t) T(t)= T(s) U(s,t).
\end{equation}
We call $\{T(t): \, t \geq 0\}$ a family of symmetries for
$\{\hat{U}(t, s): \, 0 \leq s \leq t\}$ if
for all $0 \leq s \leq t$
\[
	\hat{U}(t, s) T(s)= T(t) \hat{U}(t, s).
\]
\ed
\br
If the time-independent operator $T$ commutes with $U(s,t)$
(resp. $\hat{U}(t,s)$ for all $0 \leq s\leq t$), then
putting $T(t)=T$ for all $t\geq 0$ yields a family of symmetries.
\er
We also define symmetry on the infinitesimal level, and prove the equivalence with symmetry on the level of the evolution system.
\bd[Family of symmetries for time-dependent generators]\label{infinitesimalsym}
We call $\{T(t): \, t \geq 0\}$ a family of symmetries for $\{L_t: \, t \geq 0\}$ if $T(t)$ is differentiable in $t \geq 0$ and for all $t \geq 0$
\begin{equation}\label{eq:gen_symmetry}
	T(t) L_t - L_t T(t) = \dert T(t).
\end{equation}
We call $\{T(t): \, t \geq 0\}$ a family of symmetries for $\{\hat{L}_t: \, t\geq 0 \}$ if $T(t)$ is differentiable in $t \geq 0$ and for all $t \geq 0$
\[
	\hat{L}_t T(t) - T(t) \hat{L}_t =\dert T(t).
\]
\ed
\bp\label{equivalencesymmetries}
\begin{enumerate}
	\item [a)] $\{T(t): \, t \geq 0\}$ is a family of symmetries for  $\{L_t: \, t \geq 0\}$ if and only if it is also a family of symmetries for $\{U(s,t): \, 0 \leq s \leq t\}$ and $T(t)$ is differentiable in $t \geq 0$.
	\item [b)] $\{T(t): \, t \geq 0\}$ is a family of symmetries for $\{\hat{L}_t: \, t\geq 0 \}$ if and only if it is also a family of symmetries for $\{\hat{U}(t, s): \, 0 \leq s \leq t\}$ and $T(t)$ is differentiable in $t \geq 0$.
\end{enumerate}
\ep
\begin{proof}
	We only prove the statement for the forward process, the proof of the backward process is analogous.
	We calculate the derivative of the left hand side of \eqref{eq:evol_sys_symmetry} w.r.t. $t$, to obtain
	\[
		\pdert U(s, t) T(t) = U(s, t) L_t T(t) + U(s, t) \dert T(t).
	\]
	Setting $s = t$ gives
	\[
		\pdert U(s, t) T(t) \bigg|_{s=t} = L_t T(t) + \dert T(t).
	\]
	We also calculate the derivative of the right hand side of \eqref{eq:evol_sys_symmetry} w.r.t. $t$, namely
	\[
		\pdert T(s) U(s, t) = T(s) U(s, t) L_t.
	\]
	Choosing $s=t$ gives
	\[
		\pdert T(s) U(s, t) \bigg|_{s = t} = T(t) L_t.
	\]
	Combining these two derivatives gives that $T(t)$ is a family of symmetries for $L_t$.

	To prove that generator symmetry implies evolution system symmetry, we consider the function
	\[
		\phi(r) = U(s, s+r) T(s+r) U(s+r, t).
	\]
	\sloppy Notice that $\phi(0) = U(s, s) T(s) U(s, t) = T(s) U(s, t)$ and that $\phi(t-s) = U(s, t) T(t) U(t, t) = U(s, t) T(t)$.
	Therefore, if the derivative of $\phi$ w.r.t. $r$ is equal to 0, then we can conclude that $\phi(0) = \phi(t-s)$.
	So we calculate the derivative
	\begin{align*}
		\frac{\diff}{\diff r} \phi(r) &= U(s, s+r) \left[L_{s+r}
			T(s+r) - T(s+r) L_{s+r} + \frac{\diff}{\diff r} T(s+r) \right] U(s+r, t).
	\end{align*}
	Using \eqref{eq:gen_symmetry}, we find $\frac{\diff}{\diff r} \phi(r) = 0$, which concludes the proof.
\end{proof}
The following result shows that also in the time-inhomogeneous setting symmetries applied to duality functions yield new duality functions.
\bp \label{prop:new_dual_from_symmetry}
If $D(t)$ is a duality function between $\{U(s,t): \, 0 \leq s \leq t\}$ and $\{\hat{U}(t, s): \, 0 \leq s \leq t\}$, then we have the following:
\begin{enumerate}
	\item[a)] If $\{T(t): \, t \geq 0\}$ is a family of symmetries for $U(s,t)$, then $D(t)T(t)^{\top}$ is also a duality function between $U(s, t)$ and $\hat{U}(t, s)$, i.e., $\hat{U}(t, s)(D(s)T(s)^{\top}) = (D(t) T(t)^{\top}) U(s, t)^{\top}$.
	\item[b)] If  $\{T(t): \, t \geq 0\}$ is a family of symmetries for $\hat{U}(t, s)$, then $T(t)D(t)$ is also a duality function between $U(s, t)$ and $\hat{U}(t, s)$, i.e., $\hat{U}(t, s)(T(s) D(s)) = (T(t)D(t)) U(s, t)^{\top}$.
\end{enumerate}
\ep
\bpr
We prove the first statement, the second is proved analogously.
The statement is proven via direct calculation
\begin{align*}
	\hat{U}(t,s) (D(s) T(s)^{\top}) &= D(t) U(s, t)^{\top}
	T(s)^{\top} \\ &= D(t) (T(s) U(s, t))^{\top} \\ &= D(t) (U(s, t)T(t))^{\top} \\ &= (D(t) T(t)^{\top}) U(s, t)^{\top}.
\end{align*}
\epr

\br
Combining Proposition~\ref{equivalencesymmetries} with Proposition~\ref{prop:new_dual_from_symmetry} gives a similar statement for generators. 
\er

We have seen that it is possible to produce duality functions by acting with families of symmetries on cheap duality functions. 
It is also possible to produce a symmetry from a cheap duality function.
To make this precise, we first introduce the concept of a left and right kernel operator. 
An operator $\{T(t): \, t \geq 0\}$ acting on functions $f: S \to \R$ is in kernel form w.r.t. a measure $\{\mu_t: \, t \geq 0\}$ with right kernel $D: [0, \infty) \times S \times S$ if we can write
\[
	T(t) f(y) = \sum_{x \in S}
	D(t, y, x) f(x) \mu_t(x).
\]
Similarly we can write the operator in kernel form with left kernel $D: [0, \infty) \times S \times S$ if
\[
	T(t) f(x) = \sum_{y \in S}
	D(t, y, x) f(y) \mu_t(y).
\]

We find that the operators in kernel form given above are symmetries if and only if the transpose of the kernel used to create the kernel form is a self-duality function.
We make this statement precise in the next corollary.

\begin{corollary} \label{temp}
	Let $\{X(t): \, t \geq 0\}$ be a Markov process on the state space $S$, with evolution system $\{U(s, t): \, 0 \leq s \leq t \}$.
	Let $\{\mu_t: \, t \geq 0\}$ be a consistent trajectory of probability measures for this process.
    Let $D:[0,\infty)\times S\times S\to\R$ be a function. 
	Then the following two statements hold.
	\begin{enumerate}
		\item[a)] Define $\{T(t): \, t \geq 0\}$, with $T(t) = D(t) (D_{\mathrm{cheap}}(t))^{-1} = D(t) \text{Diag}(\mu_t)$, where $\text{Diag}(\mu_t)$ denotes a diagonal matrix with $\mu_t$ on the diagonal. 
        $T(t)$ is an operator in kernel form w.r.t. $\mu_t$ with right kernel $D(t)$. 
        $T(t)$ is a family of symmetries of $U(s, t)$ if and only if $D(t)^\top$ is a self-duality function.
		\item[b)] Define $\{T(t): \, t \geq 0\}$ with $T(t) = D(t)^\top (D_{\mathrm{cheap}}(t))^{-1} = D(t)^\top \text{Diag}(\mu_t)$.
		$T(t)$ is an operator in kernel form w.r.t. $\mu_t$ with left kernel $D(t)$. 
        $T(t)$ is a family of symmetries of $U^\rev(t, s)$ if and only if $D(t)^\top$ is a self-duality function.
	\end{enumerate}
\end{corollary}
\begin{proof}
    We only prove statement a), because the proof of the second statement is analogous.
		First we write out
		\begin{align*}
			U(s, t) T(t) f(y) &= \sum_{y', x} U(s, t)(y, y') D(t, y', x) \mu_t(x) f(x) \\
            &= \sum_x (U(s, t) D(t))(y, x) \mu_t(x) f(x).
		\end{align*}
		Second we write out, using Proposition~\ref{prop:adjoint_U}
		\begin{align*}
			T(s) U(s, t) f(y) &= \sum_{x', x} D(s, y, x') U(s, t)(x', x) \mu_s(x') f(x) \\
			&= \sum_{x', x} D(s, y, x') U^{\rev}(t, s)(x, x') \mu_t(x) f(x) \\
            &= \sum_x (D(s) U^\rev(t, s)^\top)(y, x) \mu_t(x) f(x).
		\end{align*}
    Then it is clear that if $D(t)^\top$ is a self duality function, then $T(t)$ is a family of symmetries.
    Conversely, if $T(t)$ is a family of symmetries then we find the self-duality relation $D(t)^\top U(s, t)^\top= U^{\rev}(t, s) D(s)^\top$, i.e., $D(t)^\top$ is a self duality function. 
\end{proof}

\begin{remark}
    Combining Corollary~\ref{temp} with Proposition~\ref{equivalencesymmetries} we find a similar statement for the time-dependent generators.
\end{remark}

\subsubsection{Intertwining}
An intertwining is essentially a generalization of a symmetry.
Like symmetries, intertwiners extend to the time-inhomogeneous setting.
We start by defining a time-inhomogeneous intertwiner for forward and backward evolution systems.
\begin{definition}[Intertwining family for evolution systems]
	A family of intertwiners, $\{\Lambda^{12}(t): \, t \geq 0\}$, between the evolution systems $\{U^1(s, t): \, 0 \leq s \leq t\}$ and $\{U^2(s, t): \, 0 \leq s \leq t\}$ satisfies
	\begin{equation}\label{eq:evol_sys_intertwiner}
		U^1(s, t) \Lambda^{12}(t) = \Lambda^{12}(s) U^2(s, t).
	\end{equation}
	\sloppy A family of intertwiners, $\{\Lambda^{12}(t): \, t \geq 0\}$, between the evolution systems $\{\hat{U}^1(t, s): \, 0 \leq s \leq t\}$ and $\{\hat{U}^2(t, s): \, 0 \leq s \leq t\}$ satisfies
	\begin{equation*}
		\hat{U}^1(t, s) \Lambda^{12}(s) = \Lambda^{12}(t) \hat{U}^2(t, s).
	\end{equation*}
\end{definition}
\begin{remark}
	Just as in the time-homogeneous setting, it is possible to compose intertwiners.
	So if $\Lambda^{12}(t)$ is an intertwining family between $U^1(s, t)$ and $U^2(s, t)$ and $\Lambda^{23}(t)$ is an intertwining family between $U^2(s, t)$ and $U^3(s, t)$, then $\Lambda^{13}(t) = \Lambda^{12}(t) \Lambda^{23}(t)$ is an intertwining family between $U^1(s, t)$ and $U^3(s, t)$.
	It works analogously for intertwiners between backward evolution systems.
\end{remark}
We also define an intertwiner for time-dependent generators and we show that this is equivalent to an intertwiner for an evolution system.
\begin{definition}[Intertwining family for time-dependent generators]
	A family of intertwiners, $\{\Lambda^{12}(t): \, t \geq 0\}$, between the time-dependent generators $\{L_t^1: \, t \geq 0\}$ and $\{L_t^2: \, t \geq 0\}$ is differentiable in $t \geq 0$ and satisfies
	\begin{equation}\label{eq:gen_intertwiner}
		L^1_t \Lambda^{12}(t) + \dert \Lambda^{12}(t) = \Lambda^{12}(t) L^2_t.
	\end{equation}
	A family of intertwiners, $\{\Lambda^{12}(t): \, t \geq 0\}$, between the generators $\{\hat{L}_t^1: \, t \geq 0\}$ and $\{\hat{L}_t^2: \, t \geq 0\}$ is differentiable in $t \geq 0$ and satisfies
	\begin{equation*}
		\hat{L}_t^1 \Lambda^{12}(t) = \Lambda^{12}(t) \hat{L}_t^2 + \dert \Lambda^{12}(t).
	\end{equation*}
\end{definition}
\begin{remark}
    The superscripts follow the convention of matrix indices.
    This means that the first number labels the rows and the second the columns.
    For a duality function, the row process is the backward or dual process.
    The indices compose like matrix indices and a transpose or inverse swaps them. 
\end{remark}
\begin{proposition}\label{prop:equiv_evol_sys_and_gen_intertwiner}
	\begin{enumerate}
		\item[a)] $\{\Lambda^{12}(t): \, t \geq 0 \}$ is an intertwining family between the evolution systems $\{U^1(s, t): \, 0 \leq s \leq t \}$ and $\{U^2(s, t): \, 0 \leq s \leq t \}$ and is differentiable in $t \geq 0$ if and only if it is also an intertwining family between the time-dependent generators $\{L_t^1: \, t \geq 0\}$ and $\{L_t^2: \, t \geq 0\}$.
		\item[b)] \sloppy $\{\Lambda^{12}(t): \, t \geq 0 \}$ is an intertwining family between the evolution systems $\{\hat{U}^1(t, s): \, 0 \leq s \leq t \}$ and $\{\hat{U}^2(t, s): \, 0 \leq s \leq t \}$ and is differentiable in $t \geq 0$ if and only if it is also an intertwining family between the time-dependent generators $\{\hat{L}_t^1: \, t \geq 0\}$ and $\{\hat{L}_t^2: \, t \geq 0\}$.
	\end{enumerate}
\end{proposition}
\begin{proof}
	We only prove the first statement, because the second statement is proven analogously.
	Assume that $\Lambda^{12}(t)$ is an intertwining family between $U^1(s, t)$ and $U^2(s, t)$.
	We differentiate the left hand side of \eqref{eq:evol_sys_intertwiner} w.r.t. $t$
	\begin{align*}
		\pdert U^1(s, t) \Lambda^{12}(t) &= U^1(s, t) L_t^1 \Lambda^{12}(t) + U^1(s, t) \dert \Lambda^{12}(t).
		\intertext{Setting $s=t$ gives}
		\pdert U^1(s, t) \Lambda^{12}(t) \bigg|_{s=t} &= L_t^1 \Lambda^{12}(t) + \dert \Lambda^{12}(t).
	\end{align*}
	We also calculate the derivative w.r.t. $t$ of the right hand side of \eqref{eq:evol_sys_intertwiner}
	\begin{align*}
		\pdert \Lambda^{12}(s) U^2(s, t) = \Lambda^{12}(s) U^2(s, t) L_t^2.
		\intertext{Setting $s=t$ gives}
		\pdert \Lambda^{12}(s) U^2(s, t) \bigg|_{s=t} = \Lambda^{12}(t) L_t^2.
	\end{align*}
	Combining these two derivatives gives that $\Lambda^{12}(t)$ is an intertwiner between $L_t^1$ and $L_t^2$.

	For the other direction, assume that $\Lambda^{12}(t)$ is an intertwiner between $L_t^1$ and $L_t^2$.
	We define the following function
	\[
		\phi(r) = U^1(s, s+ r) \Lambda^{12}(s+r) U^2(s+r, t).
	\]
	\sloppy Note $\phi(0) = U^1(s, s) \Lambda^{12}(s) U^2(s, t) = \Lambda^{12}(s) U^2(s, t)$ and $\phi(t-s) = U^1(s, t) \Lambda^{12}(t) U^2(t, t) = U^1(s, t) \Lambda^{12}(t)$.
	Therefore, if $\frac{\diff}{\diff r}\phi(r) = 0$, then $\phi(0) = \phi(t-s)$.
	So we calculate
	\begin{align*}
		\frac{\diff}{\diff r} \phi(r) &= U^1(s, s+r) \left[L_{s+r}^1 \Lambda^{12}(s+r) - \Lambda^{12}(s+r) L_{s+r}^2 + \frac{\diff}{\diff r} \Lambda^{12}(s+r) \right] U^2(s+r, t).
	\end{align*}
	Using \eqref{eq:gen_intertwiner}, we find $\frac{\diff}{\diff r} \phi(r) = 0$, which concludes the proof.
\end{proof}
\begin{remark}
	From Proposition~\ref{prop:equiv_evol_sys_and_gen_intertwiner} it is clear that we can also compose intertwiners for generators.
	So if $\Lambda^{12}(t)$ is an intertwining family between $L^1_t$ and $L^2_t$ and $\Lambda^{23}(t)$ is an intertwining family between $L^2_t$ and $L^3_t$, then $\Lambda^{13}(t) = \Lambda^{12}(t) \Lambda^{23}(t)$ is an intertwining family between $L^1_t$ and $L^3_t$.
	The composition of intertwining families is analogous for generators corresponding to a backward evolution system.
\end{remark}
Intertwiners provide another way of constructing duality functions.
Suppose $D(t)$ is a duality function between two processes, and suppose one of these processes is intertwined with a third process via an intertwiner.
Then this intertwiner can be composed with the duality function $D(t)$ to produce a new duality relation between the third process and the original process.
This is made precise in the proposition below.
\begin{proposition}[Duality and intertwining for evolution systems]\label{prop:duality_and_intertwining_evol_sys}
	Let $\{\hat{U}^1(t, s): \, 0 \leq s \leq t\}$ be a backward evolution system and let $\{U^2(s, t): \, 0 \leq s \leq t\}$  be a forward evolution system.
	Assume that $D^{12}(t)$ is a duality function between $U^2(s, t)$ and $\hat{U}^1(t, s)$, i.e.,
	\begin{equation*}
		D^{12}(t) U^2(s, t)^{\top} = \hat{U}^1(t, s) D^{12}(s).
	\end{equation*}
	\begin{enumerate}
		\item[a)] Let $\{\hat{U}^3(t, s): \, 0 \leq s \leq t\}$ be a backward evolution system.
		If $\Lambda^{31}(t)$ is an intertwining family between $\hat{U}^3(t, s)$ and $\hat{U}^1(t, s)$, i.e.,
		      \begin{equation*}
			      \hat{U}^3(t, s)\Lambda^{31}(s) = \Lambda^{31}(t) \hat{U}^1(t, s)
		      \end{equation*}
		      then $D^{32}(t) = \Lambda^{31}(t)D^{12}(t)$ is a duality function between $U^2(s, t)$ and $\hat{U}^3(t, s)$, i.e.,
		      \begin{equation*}
			      D^{32}(t)U^2(s, t)^{\top} = \hat{U}^3(t, s) D^{32}(s).
		      \end{equation*}
		\item[b)] Let $\{U^3(s, t): \, 0 \leq s \leq t \}$ be a forward evolution system.
		If $\Lambda^{32}(t)$ is an intertwining family between $U^3(s, t)$ and $U^2(s, t)$, i.e.,
		      \begin{equation*}
			      U^3(s, t)\Lambda^{32}(t) = \Lambda^{32}(s) U^2(s, t)
		      \end{equation*}
		      then $D^{13}(t) = D^{12}(t)\Lambda^{32}(t)^{\top}$ is a duality function between  $U^3(s, t)$ and $\hat{U}^1(t, s)$, i.e.,
		      \begin{equation*}
			      D^{13}(t)U^3(s, t)^{\top} = \hat{U}^1(t, s) D^{13}(s).
		      \end{equation*}
	\end{enumerate}
\end{proposition}
\begin{proof}
	Only the first result is proven, the second result is proved analogously.
	Using duality and the intertwining relation consecutively, we find
	\begin{equation*}
		\Lambda^{31}(t)D^{12}(t)U^2(s, t)^{\top} = \Lambda^{31}(t) \hat{U}^1(t, s)D^{12}(s) = \hat{U}^3(t, s) \Lambda^{31}(s) D^{12}(s).
	\end{equation*}
\end{proof}

Proposition~\ref{prop:duality_and_intertwining_evol_sys} shows that we can use intertwiners to find new duality functions.
Conversely, duality functions can be used to build an intertwiner.
The details can be found in the proposition below.

\begin{proposition} \label{prop:intertwiner_from_duality}
    Let $\{\hat{U}^1(t, s): \, 0 \leq s \leq t\}$ be a backward evolution system and let $\{U^2(s, t): \, 0 \leq s \leq t \}$ be a forward evolution system. 
    Let $D^{12}: [0, \infty) \times S \times S \to \R$ be a duality function between $U^2(s, t)$ and $\hat{U}^1(t, s)$, i.e.,
		      \be \label{eq:temp}
			      D^{12}(t) U^2(s, t)^\top = \hat{U}^1(t, s) D^{12}(s).
		      \ee
	\begin{enumerate}
		\item[a)]  Let $\{U^3(s, t): \, 0 \leq s \leq t \}$ be a forward evolution system.
		      And let $D^{13}: [0, \infty) \times S \times S \to \R$ be a duality function between $U^3(s, t)$ and $\hat{U}^1(t, s)$, i.e.,
		      \be \label{eq:temp_2}
		      D^{13}(t) U^3(s, t)^\top = \hat{U}^1(t, s) D^{13}(s).
		      \ee
		      \sloppy If we can also invert the matrix $D^{12}(t)$ for all $t \geq 0$, then $\Lambda^{32}(t) = D^{13}(t)^\top \left(D^{12}(t)^{-1}\right)^\top$ is an intertwiner between $U^3(s, t)$ and $U^2(s, t)$.
		\item[b)] Let $\{\hat{U}^3(t, s): \, 0 \leq s \leq t \}$ be a backward evolution system.
		      And let $D^{32}: [0, \infty) \times S \times S \to \R$ be a duality function between $U^2(s, t)$ and $\hat{U}^3(t, s)$, i.e.,
		      \[
			      D^{32}(t) U^2(s, t)^\top = \hat{U}^3(t, s) D^{32}(s).
		      \]
		      If we can also invert the matrix $D^{32}(t)$ for all $t\geq 0$, then $\Lambda^{13}(t) = D^{12}(t) D^{32}(t)^{-1}$ is an intertwiner between $\hat{U}^1(t, s)$ and $\hat{U}^3(t, s)$.
	\end{enumerate}
\end{proposition}
\begin{proof}
	We only prove the first item, because the second is proven analogously.
	Notice that we can rewrite \eqref{eq:temp} as
	\be \label{eq:temp1}
	U^2(s, t)^\top D^{12}(s)^{-1} = D^{12}(t)^{-1} \hat{U}^1(t, s).
	\ee
	We verify the definition of an intertwiner, using \eqref{eq:temp1} and \eqref{eq:temp_2}
	\begin{align*}
		U^3(s, t) \Lambda^{32}(t) &= U^3(s, t) D^{13}(t)^\top \left(D^{12}(t)^{-1}\right)^\top \\ 
        &= \left(\hat{U}^1(t, s) D^{13}(s) \right)^\top \left(D^{12}(t)^{-1}\right)^\top \\ 
        &= D^{13}(s)^\top \left(D^{12}(t)^{-1} \hat{U}^1(t, s)\right)^\top \\ 
        &= D^{13}(s)^\top \left(U^2(s, t)^\top D^{12}(s)^{-1}\right)^\top \\ 
        &= \Lambda^{32}(s) U^2(s, t).
	\end{align*}
\end{proof}

Just like it is possible to produce a symmetry from a cheap duality function, it is possible to produce an intertwiner.
We make this precise in the next corollary. 

\begin{corollary} \label{temp2}
    Let $\{U^1(s, t): \, 0 \leq s \leq t\}$ and $\{U^2(s, t): \, 0 \leq s \leq t\}$ be two forward evolution systems.
    Let $\{\mu^1_t: \, t \geq 0\}$ and $\{\mu^2_t: \, t \geq 0\}$ be consistent trajectories of probability measures for $U^1(s, t)$ and $U^2(s, t)$, respectively.
    Let $D^{12}: [0, \infty) \times S \times S \to \R$ be a function.
    Then the following two statements hold.
    \begin{enumerate}
        \item[a)] Define $\{\Lambda^{21}(t): \, t \geq 0\}$, with $\Lambda^{21}(t) = D^{12}(t)^\top \text{Diag}(\mu^1_t)$, where $\text{Diag}(\mu^1_t)$ denotes a diagonal matrix with $\mu^1_t$ on the diagonal.
        $\Lambda^{21}(t)$ is an operator in kernel form w.r.t. $\mu^1_t$ with right kernel $D^{12}(t)^\top$.
        $\Lambda^{21}(t)$ is an intertwining family between $U^2(s, t)$ and $U^1(s, t)$ if and only if $D^{12}(t)$ is a duality function between $U^2(s, t)$ and $U^{1, \rev}(t, s)$, which is the reversed evolution system with respect to $\mu^1_t$.
        \item[b)] Define $\{\Lambda^{12}(t): \, t \geq 0\}$, with $\Lambda^{12}(t) = D^{12}(t) \text{Diag}(\mu^2_t)$, where $\text{Diag}(\mu^2_t)$ denotes a diagonal matrix with $\mu^2_t$ on the diagonal.
        $\Lambda^{12}(t)$ is an operator in kernel form w.r.t. $\mu^2_t$ with right kernel $D^{12}(t)$.
        $\Lambda^{12}(t)$ is an intertwining family between $U^{1, \rev}(t, s)$ and $U^{2, \rev}(t, s)$ if and only if $D^{12}(t)$ is a duality function between $U^2(s, t)$ and $U^{1, \rev}(t, s)$. Here, $U^{1, \rev}(t, s)$ and $U^{2, \rev}(t, s)$ are reversed evolution systems with respect to $\mu^1_t$ and $\mu_t^2$, respectively. 
    \end{enumerate}
\end{corollary}
\begin{proof}
    We only prove statement a), because the proof of statement b) is analogous.
    First, we write
    \begin{align*}
        U^2(s, t) \Lambda^{21}(t) f(y)= \sum_x (U^2(s, t) D^{12}(t)^\top)(y, x) \mu_t^1(x) f(x).
    \end{align*}
    Second, we write using Proposition~\ref{prop:adjoint_U}
    \begin{align*}
        \Lambda^{21}(s) U^1(s, t) f(y) &= \sum_{x, x'} D^{12}(s, x', y) \mu_s^1(x') U^1(s, t)(x', x) f(x) \\
        &= \sum_{x, x'} D^{12}(s, x', y) U^{1, \rev}(t, s)(x, x') \mu_t^1(x) f(x) \\
        &= \sum_x (D^{12}(s)^\top U^{1, \rev}(t, s)^\top)(y, x) \mu_t^1(x) f(x).
    \end{align*}
    Then it is clear that if $D^{12}(t)$ is a duality function between $U^2(s, t)$ and $U^{1, \rev}(t, s)$, then $\Lambda^{21}(t)$ is an intertwining family.
    Conversely, if $\Lambda^{21}(t)$ is an intertwining family then we find that $D^{12}(t)$ is a duality function between $U^2(s, t)$ and $U^{1, \rev}(t, s)$.
\end{proof}

\begin{remark}
    Using Proposition~\ref{prop:equiv_evol_sys_and_gen_intertwiner} it is possible to formulate the results for evolution systems also for time-dependent generators.
\end{remark}

\subsection{Examples of symmetries of interacting particle systems}
In this section, we show that the time-inhomogeneous random walkers have time-inhomogeneous symmetries similar to the time-homogeneous case, namely the creation and annihilation operators.
We also show that there does not exist a proper time-dependent symmetry in the form of a weighted sum of creation and annihilation operators for the symmetric inclusion process (SIP) and the symmetric exclusion process (SEP). 
\subsubsection*{Example 4.1: Independent random walkers}\label{ex:time-inhomo_part_rem_op}
In this example we consider the configuration process of independent random walkers with symmetric transition rates.
We consider a (finite) vertex set $V$ and transition rates $r_t(x,y)=r_t(y,x), x,y\in V$.
Let $\{\eta(t): \, t\geq 0\}$ be the configuration process, which at each $x\in V$ counts the number of walkers $\eta_x(t)$ at time $t$.
This process has the following generator
\be\label{ltindrandwalk}
L_t=-\frac{1}{2} \sum_{x, y\in V} r_t(x, y)\left(a_y-a_x\right)\left(a_y^{\dagger}-a_x^{\dagger}\right).
\ee
Here $a_x, a_x^\dagger$ denote the particle annihilation and creation operators given by
\begin{align}\label{eq:def_part_rem_op}
	\begin{split}
		a_x f(\eta) &= \begin{cases}\eta_x f\left(\eta-\delta_x\right) & \text { if } \eta_x \geq 1 \\ 0 & \text { else, }\end{cases}\\
		a_x^{\dagger} f(\eta)&=f\left(\eta+\delta_x\right),
	\end{split}
\end{align}
where $\delta_x$ denotes one particle at position $x$.
These operators satisfy the commutation relations of the conjugate Heisenberg algebra (see e.g., \cite{redig2026duality}), i.e.,
\begin{align}\label{eq:commutation_add_rem_op}
	\begin{split}
		[a_x, a_y] &= 0, \\
		[a^{\dagger}_x, a^{\dagger}_y] &= 0, \\
		[a_x, a^{\dagger}_y] &= - \delta_{x, y}.
	\end{split}
\end{align}

It is straightforward to verify that the generator $L_t$ commutes with $\sum_x a_x$ and $\sum_x a_x^\dagger$.
We find additional time-dependent symmetries.
\bp
Define the operators
\begin{align*}
	\begin{split}
		T(t)&= \sum_{z} c_t(z) a_z, \\ T(t)^\dagger &= \sum_{z} c_t(z) a_z^\dagger,
	\end{split}
\end{align*}
with $c_t: V \to \R$ a differentiable function in $t \geq 0$.
Then, the family of operators, $\{T(t), t\geq 0\}$, is a family of symmetries for $L_t$ if and only if $c_t(z)$ satisfies the Kolmogorov backward equation, i.e.,
\begin{equation} \label{eq:der_c_t(z)}
	\pdert c_t(z)=\sum_y r_t(z, y)(c_t(z) - c_t(y)).
\end{equation}
Similarly, the family of operators  $\{T(t)^\dagger, t\geq 0\}$, is a time-dependent symmetry if and only if $c_t(z)$ satisfies the Kolmogorov forward equation, i.e.,
\[
	\pdert c_t(z)=\sum_y r_t(z, y)(c_t(y) - c_t(z)).
\]
\ep
\bpr
We prove the statement for the $T(t)$, the proof for the case of $T(t)^\dagger$ is similar.
We calculate the commutator of $T(t)$ and $L_t$
\begin{align*}
	\left[T(t), L_t \right]
	&= -\frac{1}{2} \sum_{x, y} \sum_{z} r_t(x, y) c_t(z) {\left[a_z\left(a_y-a_x\right)\left(a_y^{\dagger}-a_x^{\dagger}\right)\right.} \\
	& \qquad \left.-\left(a_y-a_x\right)\left(a_y^{\dagger}-a_x^{\dagger}\right) a_z\right] \\
	&= -\frac{1}{2} \sum_{x, y, z} r_t(x, y) c_t(z) [a_z, (a_y - a_x) (a^{\dagger}_y - a^{\dagger}_x)]. \\
	\intertext{Using that $a_x, \, a_y$ and $a_z$ commute, we find}
	\left[\sum_z c_t(z)a_z, L_t \right] &= -\frac{1}{2} \sum_{x, y, z} r_t(x, y) c_t(z) (a_y - a_x) \left([a_z, a^{\dagger}_y] - [a_z, a^{\dagger}_x] \right). \\
	\intertext{The relations given in Equation~\eqref{eq:commutation_add_rem_op} give}
	\left[\sum_z c_t(z)a_z, L_t \right] &=-\frac{1}{2} \sum_{x, y, z} r_t(x, y) c_t(z)\left(a_y-a_x\right)\left[\delta_{x, z}-\delta_{z, y}\right] \\
	&=  \sum_{y, z} r_t(z, y) c_t(z)\left(a_z-a_y\right) \\
	&= \sum_z a_z \sum_y r_t(z, y) (c_t(z) - c_t(y)). \\
\end{align*}
Using Definition~\ref{infinitesimalsym} we then find that $T(t)$ is a time-dependent symmetry if and only if we have
\[
\left[T(t), L_t \right] = \dert T(t)
\]
which leads to \eqref{eq:der_c_t(z)}.
\epr
\br
Notice that in the proof of the fact that $T(t)$ is a symmetry we only used the algebraic form of the generator, and not the particular representation of the operators.
Therefore, for any process of which the generator can be written in the form \eqref{ltindrandwalk}, we have that $T(t)$ and $T(t)^\dagger$ are time-dependent symmetries.
\er

\subsubsection*{Example 4.2: Symmetric inclusion process (SIP) and symmetric exclusion process (SEP)}\label{syminclusion}

Just as for the system of independent random walkers studied in the previous example, the generators of SIP and SEP can be written as linear combinations of operators of a Lie algebra.
For the independent random walkers, we obtained that a well-chosen time-dependent linear combination of these operators forms a proper time-dependent symmetry.
In this example, we show that such a symmetry can not exist for SIP or SEP.
We only show this result for the SIP, because the structure of the generators of both processes and the commutation relations of the operators are very similar.

We consider a finite vertex set $V$ with symmetric time-dependent transition rates $r_t(x, y) = r_t(y, x) >0$ for all $t \geq 0$ and $x, y \in V$.
Then we can write the generator of the configuration process of the SIP associated to these time-dependent edge rates as follows
(see e.g.,\ \cite{redig2026duality}, Chapter 5) with $\alpha >0$:
\[
	L_t = \frac{1}{2}\sum_{x,y} r_t(x,y)\left[K_x^+ K_y^- + K_x^- K_y^+ - 2K_x^0 K_y^0 + \frac{\alpha^2}{2}\right].
\]
The operators act on functions $f: \N \to \R$ in the following way 
\begin{align*}
	K^+ f(n) &= (\alpha + n) f(n+1), \\ K^- f(n) &= n f(n-1), \\ K^0 f(n) &= \left(\frac{\alpha}{2} + n \right) f(n).
\end{align*}
They satisfy the following commutation relations
\begin{align*}
	[K_x^\pm, K_y^0] &= \pm K_x^\pm \,\delta_{x,y} \\
	[K_x^+, K_y^-] &= 2K_x^0\,\delta_{x,y} \\
	[K_x^\pm, K_y^\pm] &= 0.
\end{align*}
Here $K_x^i$ denotes the operator $K^i$ acting
on the variable $\eta_x$.

The next proposition shows that symmetries of
the form $\sum_x\sum_{i\in \{+,-,0\}} c^i_t (x) K^i_x$ exist only when $c^i_t(x)  $ does not explicitly depend on $x$ or on $t$.
This then recovers the ``classical'' symmetries of the generator, and shows that no new time-dependent symmetries of this form arise, unlike for independent random walkers in the previous example.
\bp\label{nosymprop}
Consider operators of the form
\be\label{bopu}
	T(t) = \sum_z (c_t^+(z)\,K_z^+ + c_t^-(z)\,K_z^- + c_t^0(z)\,K_z^0).
\ee 
Then $T(t)$ is a symmetry if and only if the coefficients $c^i_t (z)$ do
not depend on time or space.
\ep 
\bpr 
To show this we calculate the commutator.
We immediately use that $K^i_x$ commutes with $K^j_y$ for $i, j \in \{+, -, 0\}$ and $x \neq y$
\begin{align}
	[L_t, T(t)] &= \frac{1}{2}\sum_{x,y} r_t(x,y)\left[K_x^+K_y^- + K_x^-K_y^+ - 2K_x^0K_y^0,\; c_t^+(x)K_x^+ + c_t^+(y)K_y^+ \right.\nonumber\\
		& \qquad \left. + c_t^-(x)K_x^- + c_t^-(y)K_y^-
		+ c_t^0(x)K_x^0 + c_t^0(y)K_y^0\right] \nonumber\\
	&= \sum_{x,y} r_t(x,y) \left[K_x^0K_y^+\left(c_t^+(y) - c_t^+(x)\right) + K_x^+K_y^0\left(c_t^+(x) - c_t^+(y)\right) \right.\nonumber\\
		& \qquad \left.+ K_x^0K_y^-\left(c_t^-(x) - c_t^-(y)\right) + K_x^-K_y^0\left(c_t^-(y) - c_t^-(x)\right) \right.\nonumber\\
		& \qquad \left. + \tfrac{1}{2}K_x^+K_y^-\left(c_t^0(x) - c_t^0(y)\right) + \tfrac{1}{2}K_x^-K_y^+\left(c_t^0(y) - c_t^0(x)\right)\right]. \label{wafa}
\end{align}
If $T(t) $ is a symmetry, then 
the expression in the right hand side of \eqref{wafa} has to equal
$- \dert (T(t))$ which is a sum of
first order terms in the Lie algebra generators 
$K^i_x, i\in\{+,-,0\}, x\in V$.
Therefore, this can only happen if all coefficients in the right hand side of \eqref{wafa} are zero, which is precisely the case when the coefficients $c^i_t (z)$ do
not depend on space.
Indeed, (only) in that case we have

\[
	[L_t, T(t)] = \dert T(t) = 0.
\]
Therefore $T(t)$ cannot explicitly depend on time.
\epr
\br 
For the SEP we have a similar representation of the generator in terms of the Lie-algebra elements
$J^+_x, J^-_x, J^0_x$ generating $\mathfrak{su}(2)$.
The same type of computation then shows that
for the SEP with time dependent symmetric edge-rates, symmetries of the form
\[
\sum_z (c_t^+(z)\, J_z^+ + c_t^-(z)\, J_z^- + c_t^0(z)\,J_z^0)
\]
exist only if $c_t^i(z)$ does not depend on space or time.
\er 

\section{Further examples}\label{sec:examples}
This section provides an example of duality for time-inhomogeneous independent random walkers, which is created by acting with a symmetry on a consistent trajectory of measures.
This construction does not work for SIP and SEP in the time-inhomogeneous setting, because as mentioned before we do not expect that there exist a proper time-dependent consistent trajectory of measures.
However, in the time-inhomogeneous setting we do still have many-to-one duality for the SIP and the SEP, where the dual process is a single random walker.
The second example shows that it is possible to obtain a consistent trajectory of measures for the SIP, by adding time-inhomogeneous reservoirs on every site.
However, this process does not have a family of symmetries.
The final example extends duality for the Ginzburg-Landau process to the time-inhomogeneous setting. 
\subsection{Independent random walkers}\label{subsec:duality_indep_walkers}
In this subsection, we show the duality relation between many independent symmetric forward random walkers and at most as many independent backward random walkers.
We consider a system of independent random walkers on the set $S$ with symmetric rates $r_t(x, y) = r_t(y, x)$ as before.
We call $\{\eta(t): \, t\geq 0\}$ the corresponding configuration process which at each $x\in S$ counts the number of walkers $\eta_x(t)$ at time $t$.
Let $\{\mu_t: \, t\geq 0\}$ be a consistent trajectory of probability measures for this process.
Let $\{\hat{\xi}(t): \, t\geq 0\}$ be the configuration process of independent backward random walkers.
We denote the process of the forward random walkers at time $t$, starting from $\eta$ at time $s$ by $\mathcal{X}(s, \eta; t)$.
We denote the reversed walkers starting from $\hat{\xi}$ at time $t$, considered at time $s$, by $\mathcal{Y}^\rev(t, \hat{\xi}; s)$.
We show that $\mathcal{X}(s, \eta; t)$ is dual to $\mathcal{Y}^\rev(t, \hat{\xi}; s)$.
We use $\mathcal{Y}^\rev$ instead of $\mathcal{X}^\rev$ because $\hat{\xi}$ is allowed to have fewer particles than $\eta$.

Below we state and prove the general result.
Then we give two concrete examples where the dual backward process consists of one or two particles.

\begin{proposition}\label{prop:duality_random_walkers}
	Let $\{\eta(t): \, t \geq 0\}$ denote the configuration of independent random walkers.
	Let $\{\mu_t: \, t\geq 0\}$ be a consistent trajectory of probability measures for one random walk.
	Let $\{\hat{\xi}(t): \, t \geq 0\}$ denote the configuration of backward independent random walkers.
	Define for $\eta \geq \hat{\xi}$
	\[
	D(t, \hat{\xi}, \eta)= \prod_{x \in S} \frac{\eta_x !}{(\eta_x-\hat{\xi}_x)! \mu_t(x)^{\hat{\xi}_x}}.
	\]
	Then we have the duality for all $ 0 \leq s \leq t$
	\[
	\E (D(t, \hat{\xi}, \mathcal{X}(s, \eta; t))) = \E (D(s, \mathcal{Y}^\rev(t, \hat{\xi}; s), \eta)).
	\]
\end{proposition}
\begin{proof}
	First notice that by the random displacement theorem (see Theorem 5.1 in \cite{last}) we have the following fact.
	If $\{\mu_t: \, t\geq 0\}$ is a consistent trajectory of probability measures, then for the configuration process $\{\eta(t): \, t\geq 0 \}$, $\prod_x \Pi(\mu_t(x))$, where $\Pi$ denotes the Poisson distribution, is a consistent trajectory of probability measures.
	From this consistent trajectory we create a cheap duality, to which we apply the exponential of the particle annihilation operator, which is a symmetry.

    So using the consistent trajectory of probability measures gives the cheap duality function

	\begin{equation}\label{cheaprw}
		D_{\mathrm{cheap}}(t, \hat{\eta}, \eta) = \prod_x D_{\mathrm{cheap}}(t, \hat{\eta}_x, \eta_x) = \prod_x \frac{\eta_x!}{\mu_t(x)^{\eta_x} e^{- \mu_t(x)}} \delta_{\eta_x, \hat{\eta}_x},
	\end{equation}
	where $\eta$ is the configuration of independent random walkers and $\hat{\eta}$ is the configuration of the backward random walkers.
	The configurations $\eta$ and $\hat{\eta}$ have the same number of particles.
	Recall the definition of the particle annihilation operator of \eqref{eq:def_part_rem_op}:
	\[a_x f(\eta)= \begin{cases}\eta_x f\left(\eta-\delta_x\right) & \text { if } \eta_x \geq 1, \\ 0 & \text { else, }\end{cases}
	\]
    where $\delta_x$ denotes one particle at position $x$.
	In this section, to alleviate notation, the operator $a$ always acts on the $\eta$ variable.
	In Example 4.1 in Subsection~\ref{ex:time-inhomo_part_rem_op} we showed that the sum of time-inhomogeneous particle annihilation operators is a symmetry of the time-inhomogeneous random walk.
	Because $\sum_x a_x$ is a symmetry, $\prod_x e^{a_x}$ is a symmetry of the time-dependent generator.
    We act with this symmetry on the $\eta$-variable of the cheap duality function
	of \eqref{cheaprw}.
    We have
	\begin{equation*}
		a^k_x D_{\mathrm{cheap}}(t, \hat{\eta}_x, \eta_x) = \frac{\eta_x!}{(\eta_x -k)!} D_{\mathrm{cheap}}(t, \hat{\eta}_x, \eta_x - k),
	\end{equation*}
	where the notation above means that the particle annihilation operator is applied $k$ times to $\eta_x$.
	As a consequence,  we find
	\begin{align*}
		\prod_x e^{a_x}
		D_{\mathrm{cheap}}(t, \hat{\eta}_x, \eta_x) &= \prod_x \sum_{k=0}^\infty \frac{1}{k!} \frac{\eta_x!}{(\eta_x - k)!} \frac{(\eta_x - k)!}{\mu_t(x)^{(\eta_x - k)} e^{-\mu_t(x)}} \delta_{\hat{\eta}_x, \eta_x - k} \\
		&= \prod_x \sum_{k=0}^\infty \frac{1}{(\eta_x - (\eta_x - k))!} \frac{\eta_x!}{\mu_t(x)^{(\eta_x - k)} e^{-\mu_t(x)}} \delta_{\hat{\eta}_x, \eta_x - k}.
	\end{align*}
	We then define $D(t, \hat{\xi}, \eta) = \prod_x e^{a_x} D_{\mathrm{cheap}}(t, \hat{\xi}_x, \eta_x)$.
    Then the previous computation yields
	\begin{equation*}
		D(t, \hat{\xi}, \eta) = \prod_x \frac{\eta_x!}{(\eta_x - \hat{\xi}_x)! \mu_t(x)^{\hat{\xi}_x} e^{-\mu_t(x)}}.
	\end{equation*}
	Note that the factor $\prod_x e^{\mu_t(x)}$ is a constant, because $\mu_t(x)$ is a probability measure and thus $\sum_x\mu_t(x)=1$ for all $t \geq 0$.
	It is always possible to multiply or divide a duality function by a constant, so we can remove this term to find the duality relation we set out to prove.    
    
    Finally, notice that the backward configuration process of independent random walkers coincides with the configuration process of independent backward random walkers.
    Indeed, if $L_t^\rev$ is the generator corresponding to the reversed process of independent random walkers then using Proposition~\ref{prop:adjoint_U} we find

    \begin{align*}
        L^\rev_t(\eta, \eta - \delta_x + \delta_y) &= \frac{\nu_t(\eta - \delta_x + \delta_y)}{\nu_t(\eta)} L_t(\eta - \delta_x + \delta_y, \eta),
    \end{align*}
    where $L_t$ is the generator corresponding to the process of independent random walkers and $\nu_t(\eta) = \prod_{z} \frac{1}{\eta_z!} \mu_t(z)^{\eta_z} e^{-\mu_t(z)}$ is the Poisson product measure corresponding to $\mu_t(x), x\in S, t\geq 0$.
    Using the definition of $\nu_t(\eta)$, we find
    \begin{align*}
        L^\rev_t(\eta, \eta - \delta_x + \delta_y) &= \frac{\frac{1}{(\eta_x-1)!} \mu_t(x)^{\eta_x -1} \frac{1}{(\eta_y+1)!}\mu_t(y)^{\eta_y+1}}{\frac{1}{\eta_x!}\mu_t(x)^{\eta_x} \frac{1}{\eta_y!} \mu_t(y)^{\eta_y}} (\eta_y+1) r_t(y, x) \\
        &= \eta_x \frac{\mu_t(y)}{\mu_t(x)} r_t(y, x) = \eta_x r^\rev_t(x, y),
    \end{align*}
    where $r_t(x, y)$ is the rate of a single random walk and $r^\rev_t(x, y)$ is the rate of a single reversed random walk.  
\end{proof}

\begin{remark}
	We illustrate the duality of Proposition~\ref{prop:duality_random_walkers} with the single and two-particle dual configurations.
	\begin{enumerate}
		\item Let $\{X^\rev(t, x; s): \, 0 \leq s \leq t \}$ denote a single random backward walker that started at time $t$ from position $x$ and we consider the walker at time $s$.
		Then we have the duality function,
		      \[
			      D(t,x,\eta)= \frac{\eta_x}{\mu_t(x)},
		      \]
		      which gives us the ``many to one duality''
		      \[
			      \E(D(t,x,\eta(t))|\eta(s)= \eta)= \E(D(s,X^{\mathrm{rev}}(t,x;s), \eta)).
		      \]
		\item Let $\{X^\rev_1(t, x_1; s): \, 0 \leq s \leq t \}$ and $\{X^\rev_2(t, x_2; s): \, 0 \leq s \leq t \}$ be two independent random backward walkers that started at time $t$ from positions $x_1$ and $x_2$, respectively, and we consider the walkers at time $s$.
		Then the duality function
		      \[D(t, x_1,x_2, \eta)=
			      \begin{cases}
				      \frac{\eta_{x_1}(t)\eta_{x_2}(t)}{\mu_t(x_1)\mu_t(x_2)} \
				      \text{if}\ x_1\not=x_2
				      \\
				      \frac{\eta_{x_1}(t)(\eta_{x_1}(t)-1)}{\mu_t(x_1)^2}\
				      \text{if} \
				      x_1=x_2,
			      \end{cases}\]
		      gives us ``many to two particles duality''
		      \[
			      \E(D(t, x_1,x_2, \eta(t))|\eta(s)=\eta)
			      =
			      \E(D(s, X_1^{\mathrm{rev}}(t,x_1;s), X_2^{\mathrm{rev}}(t,x_2; s), \eta)).
		      \]
	\end{enumerate}
\end{remark}

\subsubsection*{Many-to-one duality for the SIP and the SEP}
In the time-inhomogeneous setting, SIP and SEP are not expected to have consistent trajectories of proper time-dependent product measures
(the product of the time-independent detailed balance measure is of course a consistent trajectory).
Therefore, it is not possible to establish a duality relation similar to the one in Subsection~\ref{subsec:duality_indep_walkers} for the random walkers.
However, many-to-one duality still holds for the SIP and the SEP, which is shown via direct calculation.

\begin{proposition} \label{prop:many-to-1-SIP}
	Let $\{\eta(t): \, t \geq 0\}$ denote the configuration process which at each $x \in S$ counts the number of SIP particles, $\eta_x(t)$ at time $t$.
    The particles jump with symmetric rates $r_t(x, y) = r_t(y, x)$.
	Let $\{X(s, y; t): \, 0 \leq s \leq t\}$ be a forward random walker, with jump rate $\alpha \, r_t(x, y)$. 
    Let $\{\mu_t: \, t \geq 0\}$ be a consistent trajectory of measures for this random walker.
    Using $\mu_t$ we can construct a reversed random walk, $\{X^\rev(t, x; s): \, 0 \leq s \leq t\}$, that starts from position $x$ at time $t$, considered at time $s$.
	Then the duality function
	\begin{equation*}
		D(t, x, \eta) = \frac{\eta_x}{\mu_t(x)}
	\end{equation*}
	gives many-to-one duality
	\begin{equation*}
		\E(D(t, x, \eta(t)) | \eta(s) = \eta) = \E(D(s, X^\rev(t, x; s), \eta)).
	\end{equation*}
\end{proposition}
\begin{proof}
	The generator of the time-inhomogeneous SIP$(\alpha)$, with $\alpha>0$, is given by
	\begin{equation*}
		L_t^{\text{SIP}} f(\eta) = \sum_{x, y} r_t(x, y) \eta_x (\alpha + \eta_y) \left(f(\eta - \delta_x + \delta_y) - f(\eta)\right),
	\end{equation*}
    where $\delta_x$ denotes one particle at position $x$.
	We calculate the generator $L_t^{\text{SIP}}$ applied to the duality function $D(t, x, \eta) = \frac{\eta_x}{\mu_t(x)}$
	\begin{align*}
		L_t^{\text{SIP}} \frac{\eta_x}{\mu_t(x)} &= \sum_y \frac{r_t(x, y)}{\mu_t(x)} \left[\eta_x(\alpha + \eta_y)(\eta_x -1 - \eta_x) + \eta_y (\alpha + \eta_x)(\eta_x + 1 - \eta_x) \right]\\
		&= \sum_y \frac{r_t(x, y)}{\mu_t(x)} \alpha (\eta_y - \eta_x) \\
		&= \sum_y r_t(x, y) \alpha \left[\frac{\mu_t(y)}{\mu_t(x)} D(t, y, \eta) - D(t, x, \eta) \right].
	\end{align*}
	Note that the derivative of the duality function is
	\begin{align*}
		\pdert D(t, x, \eta) &= \pdert \frac{\eta_x}{\mu_t(x)} = -\frac{\eta_x \pdert \mu_t(x)}{\mu_t(x)^2} \\
		&= - \sum_y D(t, x, \eta) \frac{\alpha(r_t(x, y) \mu_t(y) - r_t(y, x) \mu_t(x))}{\mu_t(x)} \\
		&= - \sum_y \alpha \left[ D(t, x, \eta) r_t(x, y) \frac{\mu_t(y)}{\mu_t(x)} - D(t, x, \eta) r_t(y, x) \right].
	\end{align*}
	Combining the two calculations above gives
	\begin{align*}
		L_t^{\text{SIP}} D(t, x, \eta) + \pdert D(t, x, \eta) &= \sum_{y} r_t(x, y) \alpha \left[\frac{\mu_t(y)}{\mu_t(x)} D(t, y, \eta) - D(t, x, \eta) \right] \\
		& \qquad - \sum_y \left[ \alpha D(t, x, \eta) r_t(x, y) \frac{\mu_t(y)}{\mu_t(x)} + \alpha D(t, x, \eta) r_t(y, x) \right] \\
		\intertext{Using that $r_t(x, y) = r_t(y, x)$ provides} \\
		L_t^{\text{SIP}} D(t, x, \eta) + \pdert D(t, x, \eta) &= \sum_y \alpha r_t(x, y) \frac{\mu_t(y)}{\mu_t(x)} [D(t, y, \eta) - D(t, x, \eta)] \\
		&= \alpha \hat{L}^\text{RW}_t D(t, x, \eta),
	\end{align*}
	where $\hat{L}^\text{RW}_t f(x) = \sum_y r_t(x, y) \frac{\mu_t(y)}{\mu_t(x)} [f(y) - f(x)]$.
	By Theorem~\ref{thm:relation_gen_dual_evol_dual} generator duality implies evolution system duality and hence the result follows.
\end{proof}

\begin{remark}
    Note that the above result also holds when $\alpha$ depends on time, which can be verified with the same calculation.
\end{remark}

\begin{remark}
    We also have many-to-one duality between the SEP and one reversed random walk, namely
	\[\E(D(t, x, \eta(t)) | \eta(s) = \eta) = \E(D(s, X^\rev(t, x; s), \eta)),\]
    where $\{\eta(t): \, t \geq 0\}$ denotes the configuration process which at each $x \in S$ counts the number of SEP particles, $\eta_x(t)$ at time $t$. 
    Note that in this setting the reversed random walk has symmetric rates $r_t(x, y) = r_t(y, x)$.
	The proof is analogous to the proof for the SIP.
\end{remark}
\subsection{SIP with reservoirs}

In the time-homogeneous setting, the three particle systems IRW (independent random walkers), the SEP (symmetric exclusion) and the SIP (symmetric inclusion) have reversible product measures.
From these measures cheap dualities are obtained, which when acted upon with suitable symmetries produce the so-called triangular duality functions.
For the SEP resp. SIP with time-dependent edge rates, the reversible product measures are still the same (as in the time-homogeneous case) products of binomials resp. discrete Gamma distributions with vertex independent shape parameters.
However, for the SEP resp. SIP with time-dependent edge rates there are no proper time-dependent spatially homogeneous products of binomials resp. discrete Gamma distributions which form a consistent trajectory.
Indeed, these measures satisfy local detailed balance for every $t>0$, and therefore they cannot satisfy the Kolmogorov forward equation unless they are time independent.
In Proposition~\ref{nosymprop} we showed that
SIP with time-dependent edge rates has no proper time-dependent symmetries of the form \eqref{bopu}. We expect that more generally, there are no proper time-dependent factorized symmetries (factorized over vertices). This further suggests that there are no proper time-dependent consistent trajectories of product  measures.
\color{black}

Therefore, we consider the symmetric inclusion process with suitable time-dependent reservoirs attached to each vertex and show that this process does have a family of time-dependent product measures which is a consistent trajectory of measures.
We also prove that in that setting there is an explicit single particle duality.

We first introduce the SIP with reservoirs.
This process consists of two parts.
The first is time-homogeneous SIP dynamics on the one-dimensional lattice $\Z$  (this can be easily generalized to an arbitrary graph).
The second part is time-inhomogeneous and describes the action of time-dependent reservoirs at each vertex $x\in \Z$.
The rate at which the particles appear and disappear is determined by well-chosen time-dependent birth and death rates.
The time-dependent generator of the process is then given by
\begin{align}\label{eq:SIP_with_reservoirs}
	L_t &= L^{\text{SIP}} + \sum_{x\in \Z} L^{\text{res}, \, x}_t,
\end{align}
where for $\alpha >0$
\begin{equation*}
	\begin{split}
		L^{\text{SIP}} f(\eta) &= \sum_{x \in \Z }\left(\eta_x (\alpha + \eta_{x+1}) [f(\eta+\delta_{x+1}-\delta_x) - f(\eta)]\right. \\
		& \left. \qquad + \eta_{x+1} (\alpha + \eta_x) [f(\eta-\delta_{x+1}+\delta_x) - f(\eta)]\right)
        \color{black}
	\end{split}
\end{equation*}
describes the SIP interaction of the particles on the lattice points.
$L^{\text{res}, \, x}_t$  describes the action of the reservoir associated to vertex $x$ on the occupation variable at site $x$, and is given by
\begin{equation*}
	\begin{split}
		L^{\text{res}, \, x}_t f(\eta) &= b(t, \eta(x))[f(\eta+\delta_x) - f(\eta)] \\ & \qquad + d(t, \eta(x)) [f(\eta-\delta_x) - f(\eta)].
	\end{split}
\end{equation*}
In this process, $b(t, n) \geq 0$ and $d(t, n) \geq 0$ are the birth and death rates which are given by
\begin{align}\label{eq:time_dep_rates_birth_death}
	\begin{split}
		b(t, n) &= b(t) (\alpha + n) \\
		d(t, n) &= d(t) n.
	\end{split}
\end{align}
Here $b(t)$ and $d(t)$ are assumed to be non-negative and continuous. In the next proposition we make a precise choice such that products of discrete Gamma distributions form a consistent trajectory.
We first show that for the time-dependent birth and death process described above, a time-dependent Gamma distribution is a consistent trajectory of measures.
\begin{proposition}\label{prop:ctm_birth_death}
	For the birth and death process with rates given in \eqref{eq:time_dep_rates_birth_death},
	\begin{equation}\label{eq:c.t.m._SIP_with_reservoirs}
		\mu_t(n) = (1 - \lambda_t)^\alpha \frac{\Gamma(\alpha + n)}{\Gamma(\alpha) n!} \lambda_t^n
	\end{equation}
	is a consistent trajectory of probability measures whenever $\lambda_t \in (0, 1)$ solves
	\begin{equation}\label{eq:condition_lambda_t}
		\dert \lambda_t = (1-\lambda_t) (b(t) - d(t) \lambda_t).
	\end{equation}
\end{proposition}
\begin{proof}
	To show that \eqref{eq:c.t.m._SIP_with_reservoirs} is a consistent trajectory of measures, we show that it solves the Kolmogorov forward equation when $\lambda_t$ satisfies \eqref{eq:condition_lambda_t}.
	We consider the Kolmogorov forward equation
	\begin{equation}\label{eq:dert_mu_SIP_with_reservoirs}
		\begin{split}
			\pdert \mu_t (n) &= \mu_t(n-1) b(t) (\alpha + n -1) + \mu_t(n+1) d(t) (n+1) \\
			& \qquad - \mu_t(n) (b(t) (\alpha + n) + d(t) n).
		\end{split}
	\end{equation}
	Using \eqref{eq:c.t.m._SIP_with_reservoirs}, we find
	\begin{align*}
		\frac{\mu_t(n+1)}{\mu_t(n)} = \lambda_t \frac{\alpha + n}{n+1} && \text{and} &&
		\frac{\mu_t(n-1)}{\mu_t(n)} = \frac{n}{\lambda_t (\alpha+n-1)}.
	\end{align*}
	Using the above and \eqref{eq:dert_mu_SIP_with_reservoirs} gives
	\begin{align}\label{eq:der_mu_div_mu}
		\frac{1}{\mu_t (n)}\left(\pdert \mu_t (n)\right) &= \frac{b(t)}{\lambda_t} n + d(t) \lambda_t (\alpha + n) - b(t) (\alpha+n) - d(t) n \notag\\
		&= \alpha (d(t) \lambda_t - b(t)) + n \left(\frac{b(t)}{\lambda_t} + d(t) \lambda_t - b(t) - d(t) \right).
	\end{align}
	On the other hand, using \eqref{eq:c.t.m._SIP_with_reservoirs}, we have
	\begin{equation*}
		\log \mu_t(n) = \alpha \log(1 - \lambda_t) + n \log \lambda_t + C.
	\end{equation*}
	Therefore,
	\begin{equation*}
		\pdert \log(\mu_t(n)) = \frac{\pdert \mu_t(n)}{\mu_t(n)} = \frac{-\alpha}{1 - \lambda_t} \dert \lambda_t + \frac{n}{\lambda_t} \dert \lambda_t.
	\end{equation*}
	Comparing this with Equation~\eqref{eq:der_mu_div_mu}, we find two equations for $\dert \lambda_t$
	\begin{align*}
		\frac{1}{1- \lambda_t}\left(\dert \lambda_t\right) &= - (d(t) \lambda_t - b(t)) \\
		\frac{1}{\lambda_t}\left(\dert \lambda_t\right) &= \frac{b(t)}{\lambda_t} + d(t) \lambda_t - b(t) - d(t).
	\end{align*}
	Rewriting shows that both equations lead to the same equation for $\dert \lambda_t$ which reads
	\begin{align*}
		\dert \lambda_t &= b(t) + d(t) \lambda_t^2 - b(t) \lambda_t - d(t) \lambda_t\\
		&= (1-\lambda_t) (b(t) - d(t) \lambda_t).
	\end{align*}
\end{proof}

\begin{corollary}\label{cor:c.t.m._SIP_+_reservoirs}
	The product measure $\mu_t$ with identical marginals given in \eqref{eq:c.t.m._SIP_with_reservoirs} is a consistent trajectory of measures for the SIP with reservoirs with generator given in \eqref{eq:SIP_with_reservoirs} whenever $\lambda_t \in (0, 1)$ solves \eqref{eq:condition_lambda_t}.
\end{corollary}
\begin{proof}
	To prove this, we consider the derivative of the measure expressed in the adjoint of the generator
	\begin{align*}
		\begin{split}
			\pdert \mu_t(\eta) &= (L^{\text{SIP}})^\star \prod_{x \in \Z} \mu_t(\eta_x) + \sum_{y \in \Z}((L^{\text{res}, \, y}_t)^\star \mu_t(\eta_y)) \prod_{\substack{x \in \Z \\ x \neq y}} \mu_t(\eta_x),
		\end{split}
	\end{align*}
	$(L^{\text{SIP}})^\star \prod_{x \in \Z} \mu_t(\eta_x) = 0$ by detailed balance.
	Notice that $(L^{\text{res}, \, y}_t)^\star \mu_t(\eta_y) = (L^{\text{BD}}_t)^\star \mu_t(\eta_y)$, where $(L^{\text{BD}}_t)^\star$ is the adjoint of the generator of the birth and death process, with rates described in \eqref{eq:time_dep_rates_birth_death}.
	Combining this with Proposition~\ref{prop:ctm_birth_death} gives that $\mu_t(\eta)$ is a consistent trajectory of probability measures if $\lambda_t$ solves \eqref{eq:condition_lambda_t}.
\end{proof}
In the next proposition we obtain a form of single-particle duality for the SIP with time-dependent reservoirs defined above, i.e., the expectation of the number of particles is determined by the trajectory of a single random walk.
\bp\label{fksip}
Assume the same setting as in Proposition~\ref{prop:ctm_birth_death}.
For the process with generator \eqref{eq:SIP_with_reservoirs}, we have
\be\label{feynmansip}
\E_\eta (\eta_x(t))= e^{\int_0^t (b(r)- d(r)) \diff r} \E^{\RW}_x (\eta _{X(\alpha t)})
+ \alpha \int_0^t b(r) e^{\int_r^t (b(s)- d(s)) \diff s} \diff r
\ee
where $\E^{\RW}_x$ denotes the expectation for the nearest neighbor random walk moving at rate $1$ to each neighbor.
\ep
\bpr
Acting with the time-dependent generator on the function
$\phi(x,\eta)= \eta_x$ yields
\[
	L_t \phi(x,\eta)= \alpha \Delta_x \phi(x,\eta) + \alpha b(t) + (b(t)-d(t)) \phi(x,\eta).
\]
As a consequence, putting $\psi(t,x):=\E_\eta \eta_x(t)$ we find
\[
	\dert \psi(t,x)= \alpha \Delta_x \psi(t,x) + \alpha b(t) + (b(t)-d(t)) \psi(t,x).
\]
The solution is then given by \eqref{feynmansip} via the Feynman-Kac formula.
\epr

\subsection{Time-dependent Ginzburg-Landau process}
In the time-homogeneous setting the Ginzburg-Landau process is dual to the independent random walk, see \cite{redig2026duality}, Section 8.5.
To see this, we write the generator of the Ginzburg-Landau process as a linear combination of creation and annihilation operators.
We define two operators acting on smooth functions $f: \R \to \R$
\begin{align*}
	A_x^\dagger f(\varphi_x) &= \varphi_x f(\varphi_x) - \frac{\partial f}{\partial \varphi_x}(\varphi_x), \\ A_x f(\varphi_x) &= \frac{\partial f}{ \partial \varphi_x} (\varphi_x).
\end{align*}

These operators satisfy the canonical commutation relation $[A_x, A_y^\dagger] = I \delta_{x, y}$.
Then we can write the generator of the Ginzburg-Landau process as
\begin{equation*}
	\mathcal{L} = \sum_{x, y} \mathcal{L}_{x, y} = -\sum_{x, y} (A_x^\dagger - A_y^\dagger)(A_x - A_y).
\end{equation*}
We recall that the generator of the random walk is a linear combination of creation and annihilation operators in a discrete representation.
Therefore, we define the following operators acting on functions $f: \N \to \R$
\begin{align*}
	a_x^\dagger f(\eta_x) &= f(\eta_x + 1), \\
	a_x f(\eta_x) &= \eta_x f(\eta_x -1).
\end{align*}
These operators satisfy the dual canonical commutation relation, namely $[a_x, a_y^\dagger] = - I \delta_{x, y}$.
We write the random walk generator as a linear combination of the operators above
\begin{equation*}
	L = \sum_{x, y} L_{x, y} = -\sum_{x, y} (a_x - a_y)(a_x^\dagger - a_y^\dagger).
\end{equation*}
For the rest of this section, to alleviate notation, the operators $A$ and $A^\dagger$ always act on the $\varphi$ variable, while the operators $a$ and $a^\dagger$ act on the $\eta$ variable.
To see that the Ginzburg-Landau process is dual to the random walk, note that the operators $A_x, \, A_x^\dagger$ and $a_x, \, a_x^\dagger$ satisfy a duality relation via the Hermite orthogonal polynomials $D(\eta_x, \varphi_x) = H_{\eta_x}(\varphi_x)$, where $H_n(x) = \left(x - \derx \right)^n 1$. The duality relations are then given by
\begin{align} \label{eq:temp_4}
\begin{split}
	A_x^\dagger D(\eta_x, \varphi_x) &= a_x^\dagger D(\eta_x, \varphi_x), \\ A_x D(\eta_x, \varphi_x) &= a_x D(\eta_x, \varphi_x).
\end{split}
\end{align}

\begin{proposition}
    Using the setting described in this subsection, we have 
    \be \label{eq:temp_3}
	   \E_\varphi (D(\eta, \varphi(t)))= \hat{\E}_\eta (D(\eta(t), \varphi))
    \ee
where $\hat{\E}_\eta$
is the expectation for a system of independent {\it backward } random walkers, with time-independent duality function $D(\eta, \varphi)= \prod_x H_{\eta_x}(\varphi_x)$.
\end{proposition}
\begin{proof}
For the time-inhomogeneous Ginzburg-Landau process with generator
\begin{equation*}
	\mathcal{L}_t = \sum_{x, y} r_t(x, y) \mathcal{L}_{x, y} = -\sum_{x, y} r_t(x, y) (A_x^\dagger - A_y^\dagger)(A_x - A_y),
\end{equation*}
the operators still satisfy the duality relations in \eqref{eq:temp_4}.
Therefore we still find generator duality between the time-inhomogeneous Ginzburg-Landau process and the time-inhomogeneous random walk with generator
\begin{equation*}
	L_t = \sum_{x, y} r_t(x, y) L_{x, y} =- \sum_{x, y} r_t(x, y) (a_x - a_y)(a_x^\dagger - a_y^\dagger),
\end{equation*}
with time-independent duality function
$D(\eta, \varphi)= \prod_x H_{\eta_x}(\varphi_x)$.
From generator duality \eqref{eq:temp_3} follows. 
\end{proof}

\section*{Acknowledgment}
We used Claude Opus 5.5 to check the proofs and computations in this paper. We thank Vicente Lenz for helpful discussions.
This work was funded by the NWA ORC programme Emergence at All Scales.
\printbibliography
\end{document}